\documentclass[pdflatex,sn-mathphys-num]{sn-jnl}

\usepackage{amsmath,amssymb,amsthm}
\usepackage{mathtools}
\usepackage{tikz}
\usepackage[nameinlink,capitalise,noabbrev]{cleveref}
\crefname{thm}{Theorem}{Theorems}          \Crefname{thm}{Theorem}{Theorems}
\crefname{lem}{Lemma}{Lemmas}              \Crefname{lem}{Lemma}{Lemmas}
\crefname{prop}{Proposition}{Propositions} \Crefname{prop}{Proposition}{Propositions}
\crefname{cor}{Corollary}{Corollaries}     \Crefname{cor}{Corollary}{Corollaries}
\crefname{rem}{Remark}{Remarks}            \Crefname{rem}{Remark}{Remarks}
\crefname{defin}{Definition}{Definitions}  \Crefname{defin}{Definition}{Definitions}
\crefname{figure}{Fig.}{Figs.}             \Crefname{figure}{Fig.}{Figs.}

\theoremstyle{thmstyleone}
\newtheorem{thm}{Theorem}[section]
\newtheorem{lem}[thm]{Lemma}
\newtheorem{prop}[thm]{Proposition}
\newtheorem{cor}[thm]{Corollary}

\theoremstyle{thmstylethree}
\newtheorem{rem}[thm]{Remark}
\newtheorem{defin}[thm]{Definition}

\numberwithin{equation}{section}

\newcommand{\R}{\mathbb{R}}
\DeclareMathOperator{\dist}{dist}
\DeclareMathOperator{\sgn}{sgn}

\begin{document}

\title[Prescribed mean curvature on the ball]{Prescribed Mean Curvature on the ball via a sign change in the derivative}

\author*[1]{\fnm{\'Alvaro} \sur{Ortiz}}\email{alvaro.ortiz@uc.edu}
\author[2]{\fnm{Gonzalo} \sur{Garc\'ia}}\email{gonzalo.garcia@correounivalle.edu.co}

\affil*[1]{\orgdiv{Department of Mathematics}, \orgname{University of Cincinnati}, \orgaddress{\city{Cincinnati}, \state{OH} \postcode{45221}, \country{USA}}}
\affil[2]{\orgdiv{Departamento de Matem\'aticas}, \orgname{Universidad del Valle}, \orgaddress{\city{Cali}, \country{Colombia}}}

\abstract{We consider the problem of finding a conformal metric on the $n$-dimensional Euclidean ball $B^n$, $n\geq 3$, with zero scalar curvature in the interior and prescribed mean curvature $H$ on the boundary $\partial B^n = S^{n-1}$. For rotationally symmetric $H=H(r)$ satisfying a flatness condition of order $\alpha\in(n-2,\,n-1)$ at each of its finitely many critical points in the region where $H>0$, we prove that the sign-change of $H'(r)$ in that region is sufficient for existence. Combined with the necessary condition of Liu and Wang, this yields a characterization: in this class, the problem is solvable if and only if $H$ is not monotone where it is positive. Along the way we prove that the single blow-up point of the subcritical approximation must be a local maximum of $H$, the boundary-trace analogue of a lemma used by Chen and Li on the sphere, and, under a separation condition on the values of the local maxima, we obtain multiple solutions. Our approach adapts the variational scheme of Chen and Li for the sphere to the boundary-trace setting, using subcritical approximation and the blow-up analysis of Escobar and Garc\'ia.}



\maketitle

\section{Introduction}\label{sec:intro}

Let $(B^n, g_0)$ denote the $n$-dimensional Euclidean unit ball, $n\geq 3$, equipped with the flat metric $g_0$. The standard metric has zero scalar curvature in the interior and constant mean curvature $H_0=1$ on the boundary $\partial B^n = S^{n-1}$. A fundamental question in conformal geometry is:

\medskip
\noindent\textit{Given a function $H\colon S^{n-1}\to\R$, does there exist a metric $g=u^{4/(n-2)}g_0$ on $B^n$ with zero scalar curvature and boundary mean curvature equal to $H$?}
\medskip

\noindent This is equivalent to finding a positive function $u$ satisfying the semilinear system at the Sobolev trace critical exponent:
\begin{equation}\label{eq:main}
\begin{cases}
\Delta u = 0 & \text{in } B^n, \\[4pt]
\displaystyle\frac{\partial u}{\partial \eta}+\frac{n-2}{2}\,u = \frac{n-2}{2}\,H\,u^{n/(n-2)} & \text{on } S^{n-1}.
\end{cases}
\end{equation}

The problem was introduced by Cherrier~\cite{cherr}, who established the regularity of weak solutions, and by Escobar~\cite{con,yam,Consca,20}, who studied it as the \emph{Yamabe problem on manifolds with boundary}. Two distinct problems travel under that name. One asks for constant scalar curvature in the interior with minimal boundary; it was raised by Escobar~\cite{yam}, and most of the cases he left open were settled by Brendle and Chen~\cite{brendlechen}. The other asks for a scalar-flat metric with constant mean curvature on the boundary, which is the constant-curvature case of~\eqref{eq:main}; through the work of Escobar~\cite{Consca}, Marques~\cite{marques1,marques2}, Almaraz~\cite{almaraz} and Mayer and Ndiaye~\cite{mayerndiaye} it is now known that every compact Riemannian manifold with boundary, of dimension at least three and with finite Sobolev quotient, carries such a metric.

When the boundary curvature is prescribed rather than constant, the problem acquires a markedly different character. Escobar~\cite[Proposition~4.6]{20} showed that a Kazdan--Warner-type obstruction
\begin{equation}\label{eq:KW}
\int_{S^{n-1}} \nabla H \cdot x \; u^{2(n-1)/(n-2)}\,d\sigma = 0
\end{equation}
must hold for any positive solution $u$ of~\eqref{eq:main}, ruling out existence for large classes of functions $H$. The scalar-flat prescribed-mean-curvature problem on the ball was then studied systematically by Escobar and Garc\'ia~\cite{conformal}, who introduced the subcritical regularization scheme, characterized blow-up behavior, and proved existence results via Morse theory.

\subsection{Main results}

On the necessity side, Escudero and Garc\'ia~\cite{nota} showed, with no symmetry assumed, that if $H>0$ and $\partial H/\partial r\leq 0$ for $|x|<1$ and $H\leq 0$ for $|x|\geq 1$, then~\eqref{eq:main} has no solution. This was substantially strengthened by Liu and Wang~\cite[Theorem~1.1]{LW2022}, who proved, using fractional Laplacian techniques and the method of moving spheres, that if $H$ is continuous, radially symmetric, and monotone in the region where $H>0$, and $H\not\equiv C$, then problem~\eqref{eq:main} has no positive solution at all. Here and below, \emph{monotone in the region where $H>0$} means monotone as a function on the set $\{H>0\}$: either $H(r)\leq H(s)$ for all $r<s$ in that set, or $H(r)\geq H(s)$ for all such $r<s$. (The normalization of~\cite{LW2022} writes the boundary condition of~\eqref{eq:main} with $h\,u^{n/(n-2)}$ in place of $\gamma_n H\,u^{n/(n-2)}$, so that $h=\gamma_n H$; since $\gamma_n>0$, neither the sign nor the monotonicity hypotheses are affected.)

Our main result shows that, within a natural class of rotationally symmetric functions, this necessary condition is also sufficient.

\begin{thm}\label{thm:main}
Let $n\geq 3$ and let $H = H(r)$ be symmetric about the $x_n$-axis and positive somewhere, with $H\in C^{n-2}(S^{n-1})$ if $n\geq 4$ and $H\in C^{1,\bar\gamma}(S^{n-1})$ for some $\bar\gamma\in(0,1]$ if $n=3$. Assume that $H$ has finitely many critical points in the region where $H>0$, and suppose that near every such critical point $r_0$ the function $H$ satisfies the flatness condition
\begin{equation}\label{eq:flat}
H(r) = H(r_0) + a\,|r-r_0|^{\alpha} + k(|r-r_0|), \quad a\neq 0, \quad n-2 < \alpha < n-1,
\end{equation}
where $a=a(r_0)$, $\alpha=\alpha(r_0)$, and $k'(s) = o(s^{\alpha-1})$ as $s\to 0$; if $n\geq 5$, assume in addition that $k\in C^{n-2}$ near $0$ with
\begin{equation}\label{eq:flatplus}
k^{(s)}(t) = o\bigl(t^{\alpha-s}\bigr) \quad\text{as } t\to 0, \qquad s = 2,\ldots,n-2.
\end{equation}
If $H'(r)$ changes sign in the region where $H>0$, then problem~\eqref{eq:main} admits a positive solution $u\in C^{2,\gamma}(\overline{B^n})$ for some $\gamma\in(0,1)$, smooth in the interior of $B^n$.
\end{thm}

For $H$ of class $C^1$, the sign change of $H'$ in $\{H>0\}$ is equivalent to the failure of monotonicity there, in the sense just fixed; the verification is elementary and is carried out in \cref{rem:reading}. Together with \cref{thm:main}, the theorem of Liu and Wang therefore yields:

\begin{cor}\label{cor:char}
Let $H$ be rotationally symmetric and satisfy the regularity, finiteness, and flatness hypotheses of \cref{thm:main}, but not necessarily its sign-change hypothesis. Then the prescribed mean curvature problem~\eqref{eq:main} is solvable if and only if $H$ is not monotone in the region where $H>0$.
\end{cor}

\begin{proof}
If $H'$ changes sign in $\{H>0\}$, \cref{thm:main} produces a solution. If it does not, then $H$ is monotone in that region by \cref{rem:reading}, and $H\not\equiv C$ because a constant function has every point critical and cannot satisfy~\eqref{eq:flat}, in which $a\neq0$; so~\cite[Theorem~1.1]{LW2022} excludes every positive solution.
\end{proof}

This is the boundary-trace analogue of the characterization established by Chen and Li~\cite{pree} for the prescribed scalar curvature problem on $S^n$.

When $H$ has several positive local maxima, the minimax construction between consecutive pairs produces distinct solutions once the corresponding energy windows separate; \cref{fig:H} shows a function to which both theorems apply.

\newpage
\begin{thm}[Multiplicity]\label{thm:multi}
Let $H$ satisfy the hypotheses of \cref{thm:main} and let $r_1<r_2<\cdots<r_k$, $k\geq 2$, be the positive local maxima of $H$, ordered by latitude. Set $m_i = \min\{H(r_i),H(r_{i+1})\}$ and $\theta_n = 2^{-1/(n-2)}$. Assume that for each $i\in\{1,\ldots,k-1\}$:
\begin{enumerate}
\item no value $H(r_j)$ with $j\notin\{i,i+1\}$ lies in the interval $[\theta_n m_i,\; m_i)$;
\item the intervals $[\theta_n m_i,\; m_i]$, $i=1,\ldots,k-1$, are pairwise disjoint.
\end{enumerate}
Then problem~\eqref{eq:main} admits at least $k-1$ distinct positive solutions.
\end{thm}

The regularity, flatness and finiteness hypotheses shared by the two theorems are discussed in \cref{sec:hyp}, where the estimates that use them are available.

\begin{figure}[ht]
\centering
\begin{tikzpicture}[x=1.55cm,y=1.9cm,>=stealth,line join=round]
  \def\Hcurve{(0,-0.45) (0.35,-0.34) (0.65,-0.16) (0.85,0) (1.15,0.62) (1.45,0.92)
               (1.70,1.00) (2.00,0.92) (2.35,0.70) (2.70,0.56) (2.95,0.58)
               (3.30,0.70) (3.60,0.76) (3.90,0.70) (4.30,0.44) (4.65,0.14)
               (4.83,0) (5.20,-0.26) (5.60,-0.44) (6.00,-0.50)}
  \fill[black!8] (0,0.38) rectangle (6.15,0.76);
  \draw[->] (-0.12,0) -- (6.45,0) node[right] {$r$};
  \draw[->] (0,-0.72) -- (0,1.22) node[above] {$H$};
  \draw[very thick] (0.85,0) -- (4.83,0);
  \draw[thick,densely dotted] (0,-0.5) -- (0,-0.5);
  \draw[thick] plot[smooth,tension=0.7] coordinates \Hcurve;
  \foreach \x/\y/\n in {1.70/1.00/1, 3.60/0.76/2}{
    \fill (\x,\y) circle (1.6pt);
    \draw[densely dashed,gray] (\x,0) -- (\x,\y);
    \node[below] at (\x,-0.03) {$r_{\n}$};}
  \draw[densely dashed,gray] (0,1.00) -- (1.70,1.00);
  \draw[densely dashed,gray] (0,0.76) -- (3.60,0.76);
  \draw[densely dashed,gray] (0,0.38) -- (6.15,0.38);
  \node[left] at (-0.05,1.00) {$H(r_1)$};
  \node[left] at (-0.05,0.76) {$m_1$};
  \node[left] at (-0.05,0.38) {$\theta_n m_1$};
  \node[below left] at (-0.02,-0.03) {$0$};
  \node[below] at (6.0,-0.03) {$\pi$};
  \node at (2.85,-0.30) {$\{H>0\}$};
  \draw[<->,gray] (0.85,-0.17) -- (4.83,-0.17);
\end{tikzpicture}
\caption{A function admitted by \cref{thm:main} and \cref{thm:multi}. The derivative
$H'$ changes sign inside $\{H>0\}$, which is what \cref{thm:main} requires; the two positive
local maxima $r_1,r_2$ carry the minimax construction of \cref{sec:subcritical}, and the
shaded band is the window $[\theta_n m_1,\,m_1)$ of \cref{thm:multi}, which no other critical
value may enter. Outside $\{H>0\}$ no hypothesis is imposed, and $H$ is negative at both poles here.}
\label{fig:H}
\end{figure}

\begin{rem}[Why the two halves meet]\label{rem:reading}
For $H$ of class $C^1$, the sign change of $H'$ in $\{H>0\}$ is equivalent to the failure of monotonicity there, in the sense just fixed. In one direction, if $H'(s_2)<0$ at some $s_2\in\{H>0\}$, then $H(s)<H(s_2)$ for $s>s_2$ close to $s_2$, and such $s$ still lie in the open set $\{H>0\}$, so $H$ is not nondecreasing on it; symmetrically, $H'(s_1)>0$ prevents $H$ from being nonincreasing. In the other direction, if $\{H>0\}$ is an interval, then failure of monotonicity produces points where $H'$ takes each sign, by the mean value theorem; and if $\{H>0\}$ is disconnected, its first component $(a_1,b_1)$ has $b_1<\pi$ and $H(b_1)=0$, so $H'<0$ somewhere on it, while its last component $(a_k,b_k)$ has $a_k>0$ and $H(a_k)=0$, so $H'>0$ somewhere on it. In particular monotonicity on $\{H>0\}$ forces that set to be connected and to contain a pole, which is exactly the configuration reduced to in~\cite[Section~2]{LW2022}. Finally, the proviso $H\not\equiv C$ of~\cite{LW2022} is automatic in the class of \cref{cor:char}: a constant function has every point critical and cannot satisfy~\eqref{eq:flat}, in which $a\neq 0$.
\end{rem}

\subsection{Existing results and the role of flatness conditions}

A rich body of work has developed around equation~\eqref{eq:main} under various \emph{non-degeneracy} and \emph{flatness} conditions on $H$. The predominant approach, originating in Bahri's theory of critical points at infinity, was pioneered for problem~\eqref{eq:main} by Abdelhedi, Chtioui, and Ould Ahmedou~\cite{ACO2008,ACO2009}, who used Morse-theoretic methods to establish existence and multiplicity results in dimensions $n=4$ and $n\geq 3$, respectively. In~\cite{AC2010}, Abdelhedi and Chtioui proved that under the flatness condition
\begin{equation}\label{eq:fb}
H(x) = H(y) + \sum_{i=1}^{n-1} b_i |x_i - y_i|^{\beta} + R(x-y), \quad b_i \neq 0, \quad \beta \in (n-2,\, n-1), \tag{$f_\beta$}
\end{equation}
at each critical point $y$ of $H$ (with $\sum_{i} b_i\neq 0$ and appropriate remainder), an Euler--Hopf index-counting formula gives both existence and generic multiplicity estimates.

This line of investigation was extended to the complementary flatness window $\beta\in(1,\,n-2]$ by Bensouf~\cite{bensouf}, to the endpoint $\beta = n-2$ by Al-Ghamdi, Chtioui, and Sharaf~\cite{ACS2013}, and beyond $\beta = n-1$ by Sharaf~\cite{sharaf}. More recently, Fourti~\cite{fourti} obtained existence results under pinching conditions, while Djadli, Malchiodi, and Ould Ahmedou~\cite{DMO2004} treated the four-dimensional ball $B^4$ using topological methods.

Via different methods, Chang, Xu, and Yang~\cite{chang} established perturbative existence using degree theory when $H$ is close to a constant, and Xu and Zhang~\cite{XZ2016} extended this to a quantitative closeness criterion using the mean curvature flow. In the presence of symmetry, Ho~\cite{ho} proved existence via a flow method under the assumption that $H>0$ everywhere, $H$ satisfies a pinching condition, and the Laplacian $\Delta H$ is positive at the fixed-set maximum. In an unpublished preprint, Zhang~\cite{zhang} treats sign-changing prescribed functions by a negative gradient flow, under a pinching condition on $\max|H|$ against the mean value of $H$ together with either a Morse index-counting hypothesis or an invariance hypothesis. Wang and Zhao~\cite{WZ2013} constructed infinitely many non-rotationally symmetric solutions via Lyapunov--Schmidt reduction when $H$ has a local maximum satisfying certain expansion conditions.

In a parallel and rapidly developing direction, Cruz-Bl\'azquez, Malchiodi, and Ruiz~\cite{CMR2022} initiated the study of prescribing \emph{both} scalar curvature $K\neq 0$ and boundary mean curvature simultaneously, with important subsequent contributions by Battaglia, Pu, and Vaira~\cite{BPV2025}. These results primarily address the case $K<0$; the scalar-flat case $K=0$ of the present paper sits at the boundary of their framework.

The compactness and blow-up analysis underlying all of these approaches was developed, for the constant-mean-curvature problem, by Felli and Ould Ahmedou~\cite{FO2003}, who established $C^2$ compactness and computed the total Leray--Schauder degree on locally conformally flat manifolds with umbilic boundary.

\subsection{Comparison with existing results}

The hypotheses of \cref{thm:main} are essentially disjoint from all of the results cited above, in two respects.

First, condition~\eqref{eq:fb} requires the coefficients $b_i$ to be nonzero in \emph{every} coordinate direction, forcing each critical point of $H$ to be \emph{isolated}. However, a rotationally symmetric function $H = H(r)$ on $S^{n-1}$ with a critical latitude $r = r_0 \in (0,\pi)$ possesses an entire $(n-2)$-dimensional sphere of critical points along which all tangential derivatives vanish identically, so condition~\eqref{eq:fb} structurally fails and the Euler--Hopf index-counting framework of~\cite{AC2010,bensouf,sharaf,ACS2013,fourti} does not apply. This is precisely the gap identified by Chen and Li~\cite{pree} on $S^n$ for the prescribed scalar curvature problem: rotationally symmetric functions have degenerate critical manifolds that fall outside the purview of Y.\,Y.~Li's isolated-critical-point theory~\cite{Li}. \cref{thm:main} fills this gap for the boundary problem on the ball.

Second, unlike all existing results via the $(f_\beta)$ framework, which require $H>0$ everywhere on $S^{n-1}$, and unlike Ho's flow result~\cite{ho}, which requires $H>0$ together with a global pinching bound and $\Delta H>0$ at the fixed-set maximum, our hypotheses impose no sign restriction on $H$ outside the region where it is positive. A function satisfying our flatness~\eqref{eq:flat} with $\alpha>2$ (which holds for $n\geq 5$) has vanishing Laplacian at the critical point and thus cannot satisfy Ho's condition. The preprint~\cite{zhang} also allows sign change, but every one of its results requires the two global conditions $\int_{S^{n-1}}H\,d\sigma>0$ and $\max_{S^{n-1}}|H| < 2^{1/(n-1)}\,\overline{H}$, where $\overline{H} = |S^{n-1}|^{-1}\int_{S^{n-1}}H\,d\sigma$ is the mean value of $H$; our theorem needs neither. This includes the results that assume invariance under a group of isometries, which are the ones closest to our setting. (The exponent is stated here in our dimension convention; that preprint works on $B^{n+1}$ with boundary $S^{n}$ and writes $2^{1/n}$.)

\subsection{Strategy of proof}

The proof adapts the approach of Chen and Li~\cite{pree} for the sphere to the boundary-trace setting, using the analytic framework of Escobar and Garc\'ia~\cite{conformal}. The argument has two stages.

In \cref{sec:subcritical}, we solve the \emph{subcritical} problem ($1<p<\tau:= n/(n-2)$) by a variational argument with two cases, following the scheme of~\cite{pree}. If $H$ has at least two positive local maxima, the functional $J_p(u) = \int_{S^{n-1}} H\,u^{p+1}\,d\sigma$, restricted to the energy sphere $S$, admits two ``concentration neighborhoods'' $\Sigma_1$, $\Sigma_2$ centered at the two positive local maxima of least value. When $H'(r)$ changes sign, the functional values on $\partial\Sigma_i$ are strictly below the concentration levels $H(r_i)|S^{n-1}|$, and the mountain-pass theorem yields a critical point $u_p$ with energy below both levels; a two-bubble path (\cref{lem:lower}) bounds the minimax level away from zero, which makes the Lagrange multiplier positive. If $H$ has exactly one positive local maximum, the sign change forces it to lie at an interior latitude, with each pole either non-positive or a positive local minimum (\cref{lem:config}); in this case we maximize $J_p$ over the rotationally symmetric functions in $S$ and invoke the principle of symmetric criticality~\cite{palais}.

In \cref{sec:apriori}, we establish a priori bounds on the subcritical solutions $\{u_p\}$ as $p\to\tau$, partitioning $S^{n-1}$ into three regions: where $H$ is negative (a Harnack inequality obtained by the Kelvin-transform method of Chen and Li~\cite{chen}), where $H$ is small (rescaling at the natural boundary scale and Liouville-type arguments), and where $H$ is positive. In the positive region, the blow-up theory of Escobar and Garc\'ia~\cite{conformal} shows that any blow-up consists of a single simple bubble (\cref{lem:EG}). We then prove that the blow-up point must be a critical point and indeed a \emph{local maximum} of $H$ (\cref{lem:location}); this is the boundary-trace analogue of a lemma used by Chen and Li~\cite{pree} on the sphere, where it is quoted from Y.\,Y.~Li~\cite{Li}, and we give a self-contained proof from a global Pohozaev identity on the half-space. In the mountain-pass case the energy bound~\eqref{eq:cpbound} then rules out even this single blow-up; in the rotationally symmetric case the blow-up point would have to be a pole, which is not a local maximum. The Arzel\`a--Ascoli theorem then provides convergence to a solution of~\eqref{eq:main}.

The flatness window $n-2 < \alpha < n-1$ is determined by two analytic constraints: the lower bound $\alpha > n-2$ is what the simple-blow-up theory of~\cite{conformal} requires (Theorems~4.8, 4.10 and~4.12 there), and the upper bound $\alpha < n-1$ makes the concentration term $-C\lambda^\alpha$ dominate the $O(\lambda^{n-1})$ far-field remainders, both in the expansion of $J_p$ near a bubble (\cref{jp}) and in the Pohozaev balance of \cref{lem:location}. This is the same window that appears in~\cite{AC2010}; the coincidence reflects the common analytic origin in the bubble estimates, not a dependence between the two approaches.

\section{Preliminaries}\label{sec:prelim}

Let $\gamma_n = \frac{n-2}{2}$ and $\tau = \frac{n}{n-2}$. We consider the functionals
\[
J_p(u) = \int_{S^{n-1}} H\,(u_+)^{p+1}\,d\sigma, \qquad
E(u) = \int_{B^n} |\nabla u|^2\,dv + \gamma_n\int_{S^{n-1}} u^2\,d\sigma,
\]
where $u_+ = \max\{u,0\}$, and the constraint set
\[
S = \bigl\{u\in H^1(B^n) : E(u) = \gamma_n|S^{n-1}|\bigr\}.
\]
Since every function appearing in the constructions below is nonnegative, we write $u^{p+1}$ for $(u_+)^{p+1}$ whenever $u\geq0$; the positive part matters only in \cref{lem:multiplier}, where it makes critical points of $J_p$ on $S$ automatically nonnegative. A scalar multiple of a critical point of $J_p$ on $S$ solves the subcritical version of~\eqref{eq:main}:
\begin{equation}\label{eq:subcrit}
\begin{cases}
\Delta u = 0 & \text{in } B^n, \\[4pt]
\displaystyle\frac{\partial u}{\partial \eta}+\gamma_n\,u = \gamma_n\,H\,u^p & \text{on } S^{n-1},
\end{cases}
\end{equation}
for $1<p\leq\tau$.

We fix coordinates so that the south pole of $S^{n-1}$ is at the origin $\mathcal{O}$ and the center of $B^n$ is at $(0,\ldots,0,1)$. The center of mass of $u$ is
\[
q = q(u) = \frac{\int_{B^n} z\,u^{2\tau}(z)\,dv}{\int_{B^n} u^{2\tau}(z)\,dv}.
\]
For $q\in B^n$, let $u_q$ denote the standard solution of~\eqref{eq:subcrit} with $H\equiv 1$ and center of mass at $q$; that is, $u_q$ satisfies
\begin{equation}\label{eq:standard}
\begin{cases}
\Delta u_q = 0 & \text{in } B^n, \\[4pt]
\displaystyle\frac{\partial u_q}{\partial \eta}+\gamma_n\,u_q = \gamma_n\,u_q^\tau & \text{on } S^{n-1}.
\end{cases}
\end{equation}
Write $\tilde{q}$ for the intersection of $S^{n-1}$ with the ray from the center through $q$. The solutions $u_q$ depend on two parameters: the point $\tilde{q}$ and a number $\beta$ with $1\leq\beta<\infty$. When $\tilde{q}= 0$, the solutions are given by
\[
u_q(z) = \left(\frac{4\beta}{(\beta-1)^2\|z\|^2+4z_n(\beta-1)+4}\right)^{(n-2)/2},
\]
where $z_n$ is the last component of $z$. These solve the problem of prescribing zero scalar curvature and mean curvature $H\equiv 1$ on $B^n$. For $z\in S^{n-1}$, this simplifies to
\[
u_q(z) = \left(\frac{\beta}{\frac{(\beta^2-1)\|z\|^2}{4}+1}\right)^{(n-2)/2}.
\]
In spherical coordinates $(r,\theta)$ on $S^{n-1}$ centered at the south pole, with $0\leq r\leq\pi$ and $\theta\in S^{n-2}$, we may write
\[
u_q(z) = \left(\frac{\beta}{(\beta^2-1)\sin^2\frac{r}{2}+1}\right)^{(n-2)/2}.
\]
Setting $\lambda = 1/\beta\in(0,1]$, this becomes
\begin{equation}\label{eq:uq}
u_q(z) = \left(\frac{\lambda}{\lambda^2\cos^2\frac{r}{2}+\sin^2\frac{r}{2}}\right)^{(n-2)/2}.
\end{equation}
Using the volume relation $|S^{n-1}| = |S^{n-2}|\int_0^\pi \sin^{n-2}(t)\,dt$, one verifies that $u_q\in S$.

The conformal dilation $\varphi = \psi^{-1}\circ T\circ\psi$, where $\psi(z) = 4(z-N)/|z-N|^2+N$ and $T(x,t) = (\beta x,\, \beta(t+1)-1)$, maps $B^n$ to itself and induces the isometry $T_\varphi(u) = u\circ\varphi\cdot[\det(d\varphi)]^{(n-2)/(2(n-1))}$ on $H^1(B^n)$: it preserves $E(\cdot)$ and the critical trace integral $\int_{S^{n-1}}u^{2(n-1)/(n-2)}\,d\sigma$, the exponent $(n-2)/(2(n-1))$ being exactly the weight that makes the latter invariant (see~\cite{con}). It does not preserve $J_p$ for $p<\tau$, and only the $E$-isometry is used below.

\begin{defin}\label{def:sigma}
For $\rho_0>0$ sufficiently small, define
\[
\Sigma = \bigl\{u\in S : |q(u)|\leq\rho_0,\; \|v\| = \min_{t,q}\|u-tu_q\|\leq\rho_0,\; t\in\R\bigr\},
\]
the set of functions in $S$ whose center of mass is near the south pole $\mathcal{O}$.
\end{defin}

The following lemmas, whose proofs mirror those given in~\cite{pree} for $S^n$, are used in the estimates of \cref{sec:subcritical}.

\newpage
\begin{lem}[Center of mass]\label{lem:cmasa}
For $|q|$ sufficiently small:
\begin{enumerate}
\item $|q|^2 \leq C(|\tilde{q}|^2 + \lambda^4)$;
\item $\rho_0 \leq |q| + C\|v\|$.
\end{enumerate}
\end{lem}

\begin{lem}[Orthogonality]\label{lem:orto}
If $u\in\Sigma$ and $v = u - t_0 u_q$ as in \cref{def:sigma}, then $u_q$ and $v$ are orthogonal with respect to $E(\cdot)$, and
\[
\int_{S^{n-1}} u_q^\tau\,v\,d\sigma = 0.
\]
\end{lem}

\begin{lem}\label{lem:orto2}
If $u\in\Sigma$ and $v = u - t_0 u_q$, then $T_\varphi v$ is orthogonal to the constants and to the coordinate functions with respect to the inner product $E(\cdot,\cdot)$:
\[
\langle T_\varphi v,\, 1\rangle_{E} = 0 \quad\text{and}\quad \langle T_\varphi v,\, x_i\rangle_{E} = 0, \qquad 1\leq i\leq n.
\]
\end{lem}

We fix the following notation for the rest of the paper. For $x\in\partial\R^n_+$ and $r>0$ we
write $B_r^+(x) = B_r(x)\cap\R^n_+$ for the open half-ball, $\partial'B_r^+(x) = \partial
B_r^+(x)\cap\partial\R^n_+$ for its flat part, and $\partial''B_r^+(x) = \partial B_r^+(x)\cap\R^n_+$
for its spherical part; the flat part is denoted $\partial^0$ in~\cite{conformal}. Interiors and
boundaries of subsets of $S$, such as $\mathring{\Sigma}$ and $\partial\Sigma$, are taken relative
to $S$; by \cref{def:sigma}, $u\in\partial\Sigma$ forces $|q(u)|=\rho_0$ or
$\|v\|=\rho_0$. We write $\Sigma(r_0)$ for the set of \cref{def:sigma} with the south pole
replaced by the latitude $r_0$, and abbreviate $\Sigma_i=\Sigma(r_i)$. Finally, the standard
solutions are indexed in two ways: as $u_q$, by the center of mass $q$, and as
$u_{\lambda,\tilde{q}}$, by the concentration parameter $\lambda$ of~\eqref{eq:uq} together with the
point $\tilde{q}\in S^{n-1}$. The two notations denote the same family.

\section{Subcritical solutions via the mountain pass theorem}\label{sec:subcritical}

In this section we prove existence of a solution to~\eqref{eq:subcrit} for each $1<p<\tau$. The estimates of \cref{sec:estimates} serve the mountain-pass construction used when $H$ has at least two positive local maxima; the case of a single positive local maximum is treated directly in \cref{sec:existence}. We present the estimates near the south pole $\mathcal{O}$, where we assume $H$ has a positive local maximum; the estimates near any other positive local maximum are identical.

By the flatness hypothesis~\eqref{eq:flat}, near $\mathcal{O}$ we have $H(r) = H(0) - a\,r^{\alpha} + k(r)$ for some $a>0$ and $n-2<\alpha<n-1$, with $k$ as in~\eqref{eq:flat}--\eqref{eq:flatplus}; the coefficient is negative because $\mathcal{O}$ is a local maximum.

\subsection{Estimates on \texorpdfstring{$J_p$}{J\_p}}\label{sec:estimates}

\begin{prop}\label{prop:sup}
For all $\delta_1>0$, there exists $p_1\leq\tau$ such that for all $p_1\leq p\leq\tau$,
\[
\sup_{\mathring{\Sigma}} J_p(u) > H(0)|S^{n-1}| - \delta_1.
\]
\end{prop}

\begin{proof}
We show that $J_\tau(u_{\lambda,\mathcal{O}})\to H(0)|S^{n-1}|$ as $\lambda\to 0$. Fix a small $\theta_0 > 0$ and write
\begin{align*}
J_\tau(u_{\lambda,\mathcal{O}}) &= |S^{n-2}|\int_0^\pi \frac{\lambda^{n-1}H(r)\sin^{n-2}r}{(\lambda^2\cos^2\frac{r}{2}+\sin^2\frac{r}{2})^{n-1}}\,dr \\
&= H(0)|S^{n-1}| + |S^{n-2}|\bigl\{-I_1 + I_2\bigr\},
\end{align*}
where
\begin{gather*}
I_1 = \left|\int_0^{\theta_0} \frac{a\,r^\alpha\,\lambda^{n-1}\sin^{n-2}r}{(\lambda^2\cos^2\frac{r}{2}+\sin^2\frac{r}{2})^{n-1}}\,dr\right|, \\
I_2 = \left|\int_{\theta_0}^\pi \frac{\lambda^{n-1}(H(r)-H(0))\sin^{n-2}r}{(\lambda^2\cos^2\frac{r}{2}+\sin^2\frac{r}{2})^{n-1}}\,dr\right|.
\end{gather*}
For $I_1$ we bound the denominator by splitting its two summands rather than discarding one of them: for $A,B>0$,
\[
(A+B)^{n-1}\;\geq\;A^{(n-2)/2}B^{n/2}, \qquad \tfrac{n-2}{2}+\tfrac{n}{2}=n-1.
\]
Taking $A=\lambda^2\cos^2\frac{r}{2}$ and $B=\sin^2\frac{r}{2}$ leaves a single power of $\lambda$ together with a factor $\cos^{n-2}\frac{r}{2}$, which cancels against $\sin^{n-2}r = 2^{n-2}\sin^{n-2}\frac{r}{2}\cos^{n-2}\frac{r}{2}$. Hence
\begin{align*}
I_1 &\leq c\,\lambda \int_0^{\theta_0} \frac{r^{\alpha}\,(\sin\frac{r}{2})^{n-2}}{(\sin^2\frac{r}{2})^{n/2}}\,dr
= c\,\lambda \int_0^{\theta_0} \frac{r^{\alpha}}{\sin^2\frac{r}{2}}\,dr
\leq c_1\,\lambda \int_0^{\theta_0} r^{\alpha-2}\,dr = c_1\,\lambda\,\frac{\theta_0^{\alpha-1}}{\alpha-1},
\end{align*}
using $\sin\frac{r}{2}\geq\frac{r}{\pi}$ on $[0,\pi]$ in the last inequality; this converges to $0$ since $\alpha>1$. The cruder bound $(\lambda^2\cos^2\frac{r}{2}+\sin^2\frac{r}{2})^{n-1}\geq\sin^{2n-2}\frac{r}{2}$ would instead leave $\int_0^{\theta_0}r^{\alpha-n}\,dr$, which diverges precisely because $\alpha<n-1$. A similar computation gives $I_2\leq C_2\,\lambda\,(\pi-\theta_0)\to 0$.

Choosing $\lambda_0$ so that $u_{\lambda_0,\mathcal{O}}\in\mathring{\Sigma}$ and $J_\tau(u_{\lambda_0,\mathcal{O}})>H(0)|S^{n-1}|-\delta_1/2$, the continuity of $J_p$ in $p$ yields the result.
\end{proof}

\begin{lem}\label{jp}
For $n-2 < \alpha < n-1$, $p$ sufficiently close to $\tau$, and $\lambda$, $|\tilde{q}|$ sufficiently small,
\[
J_p(u_{\lambda,\tilde{q}}) \leq \bigl(H(0)-C_1|\tilde{q}|^\alpha\bigr)|S^{n-1}|(1+o_p(1)) - C_1\lambda^{\alpha+\delta_p},
\]
where $\delta_p = \tfrac{(n-2)(\tau-p)}{2}>0$ and $o_p(1)\to 0$ as $p\to\tau$.
\end{lem}

\begin{proof}
Write $H = H(0)-a\,r^\alpha+k(r)$ near $\mathcal{O}$, $r=\dist(\cdot,\mathcal{O})$, and split $J_p(u_{\lambda,\tilde{q}})$ over $\{r\leq\sigma_0\}$ and $\{r\geq\sigma_0\}$. The substitution $t=\tan\frac{\rho}{2}$, $t=\lambda s$ in geodesic polar coordinates centered at $\tilde q$ yields the subcritical mass and moment laws
\[
\int_{S^{n-1}}u_{\lambda,\tilde{q}}^{\,p+1}\,d\sigma = |S^{n-1}|\,\lambda^{\delta_p}\bigl(1+O(\tau-p)\bigr),
\qquad
\int_{S^{n-1}}r^{\alpha}\,u_{\lambda,\tilde{q}}^{\,p+1}\,d\sigma \asymp \lambda^{\delta_p}\bigl(\lambda^{\alpha}+|\tilde{q}|^{\alpha}\bigr),
\]
uniformly for $p$ near $\tau$, the convergence of the transversal tail requiring exactly $\alpha<n-1$. Together with $|k(r)|=o(r^{\alpha})$ and the far-field bound $\int_{\{r\geq\sigma_0\}}|H|u^{p+1}\,d\sigma = \lambda^{\delta_p}O(\lambda^{\alpha+\kappa})$, $\kappa=\tfrac{n-1-\alpha}{2}$, this gives
\[
J_p(u_{\lambda,\tilde{q}}) \leq \lambda^{\delta_p}\Bigl[\bigl(H(0)-C_1|\tilde{q}|^{\alpha}\bigr)|S^{n-1}|\bigl(1+o_p(1)\bigr)-C_1\lambda^{\alpha}+C_2\lambda^{\alpha+\kappa}\Bigr].
\]
Since the bracket is positive for $\sigma_0$ small and $\lambda^{\delta_p}\leq1$, the stated inequality follows; the passage from $\lambda^{\alpha}$ inside the bracket to the weaker $\lambda^{\alpha+\delta_p}$ in the statement is harmless, since $\lambda\leq 1$ and $\delta_p\geq 0$. The complete computation, with the two-sided moment bounds, parallels~\cite[Lemma~2.1]{pree} on $S^n$, with the boundary-trace exponents in place of the interior ones; the form of the conclusion, with $\lambda^{\alpha+\delta_p}$ rather than $\lambda^\alpha$, is the one used there.
\end{proof}

\begin{rem}\label{rem:deltap}
The exponent $\delta_p=\tfrac{(n-2)(\tau-p)}{2}$ satisfies $\delta_p>\tau-p$ only for $n\geq5$; for $n=3,4$ it is positive but smaller than $\tau-p$. Only $\delta_p\geq0$ is used in the sequel.
\end{rem}

\begin{prop}\label{prop:bdy}
There exist positive constants $\rho_0$, $p_0$, $\delta_0$ such that for all $p\geq p_0$ and $u\in\partial\Sigma$,
\[
J_p(u) \leq H(0)|S^{n-1}| - \delta_0.
\]
\end{prop}

\begin{proof}
Define $\bar{H}(x) = H(x)$ in $B_{2\rho_0}(0)$ and $\bar{H}(x) = m_0 := H|_{\partial B_{2\rho_0}(0)}$ outside, and set $\bar{J}_p(u) = \int_{S^{n-1}}\bar{H}\,u^{p+1}\,d\sigma$.

\textbf{Step~1.} The difference $|J_p(u)-\bar{J}_\tau(u)|$ is controlled by $|\bar{J}_p(u)-\bar{J}_\tau(u)| + |J_p(u)-\bar{J}_p(u)|$. For the first term, H\"older's inequality with exponents $\tfrac{\tau+1}{p+1}$ and $\tfrac{\tau+1}{\tau-p}$, together with $\bar{H}^{(\tau+1)/(p+1)}\leq\bar{H}\,(\sup\bar{H})^{(\tau-p)/(p+1)}$, gives
\[
\bar{J}_p(u) \;\leq\; \bar{J}_\tau(u)\cdot\bigl(\bar{J}_\tau(u)\bigr)^{-\frac{\tau-p}{\tau+1}}\bigl((\sup\bar{H})\,|S^{n-1}|\bigr)^{\frac{\tau-p}{\tau+1}}.
\]
On $\partial\Sigma$ the factor $\bar{J}_\tau(u)$ is bounded above by $(\sup\bar{H})|S^{n-1}|$ and below by a positive constant, since $\bar{H}\geq c>0$ on $S^{n-1}$ for $\rho_0$ small and $\int_{S^{n-1}}u^{\tau+1}\,d\sigma\geq c'$ near the bubble family; both exceptional factors are therefore of the form $A^{\tau-p}$ with $A$ in a fixed compact subset of $(0,\infty)$, hence $1+o_p(1)$ uniformly. Thus
\begin{equation}\label{eq:Jbar-sub}
\bar{J}_p(u) \leq \bar{J}_\tau(u)(1+o_p(1)),
\end{equation}
where $o_p(1)\to 0$ as $p\to\tau$. For the second term:
\begin{equation}\label{eq:Jp-Jbar}
|J_p(u)-\bar{J}_p(u)| = \int_{S^{n-1}\setminus B_{2\rho_0}(0)}|H-m_0|\,u^{p+1}\,d\sigma \leq C_1\int_{S^{n-1}\setminus B_{2\rho_0}(0)} u^{p+1}\,d\sigma.
\end{equation}
Writing $u = t_0 u_q + v$, we bound
\[
C_1\int_{S^{n-1}\setminus B_{2\rho_0}(0)} u^{p+1}\,d\sigma \leq C_2\int_{S^{n-1}\setminus B_{2\rho_0}(0)} (t_0 u_q)^{p+1}\,d\sigma + C_2\int_{S^{n-1}\setminus B_{2\rho_0}(0)} v^{p+1}\,d\sigma.
\]
A direct computation gives the first integral $\leq C_3\lambda^{(n-1)-\delta_p}$. For the second, we apply the Beckner--Escobar Sobolev trace inequality:
\[
\left(\int_{S^{n-1}} v^{p+1}\,d\sigma\right)^{1/(p+1)} \leq C\left(\int_{B^n}|\nabla v|^2\,dv + \gamma_n\int_{S^{n-1}} v^2\,d\sigma\right)^{1/2} = C\|v\|.
\]
Hence
\[
|J_p(u)-\bar{J}_p(u)| \leq C_3\lambda^{(n-1)-\delta_p} + C\|v\|^{p+1}.
\]
Since both $\lambda^{(n-1)-\delta_p}$ and $\|v\|^{p+1}$ are small, using~\eqref{eq:Jbar-sub} and~\eqref{eq:Jp-Jbar}, the difference between $J_p(u)$ and $\bar{J}_\tau(u)$ is controlled.

\textbf{Step~2.} We estimate $\bar{J}_\tau(u)$. Let $u = v + t_0 u_q\in\partial\Sigma$. By \cref{lem:orto}, $v$ and $u_q$ are orthogonal with respect to $E(\cdot)$; that is,
\[
0 = \int_{B^n} \nabla(u-t_0 u_q)\cdot\nabla u_q\,dv + \gamma_n\int_{S^{n-1}}(u-t_0 u_q)\,u_q\,d\sigma,
\]
so $t_0 E(u_q) = \int_{B^n}\nabla u\cdot\nabla u_q\,dv + \gamma_n\int_{S^{n-1}} u\,u_q\,d\sigma$. Hence
\begin{align*}
\|v\|^2 &= E(u-t_0 u_q) = E(u) + t_0^2 E(u_q) - 2t_0\left(\int_{B^n}\nabla u\cdot\nabla u_q\,dv + \gamma_n\int_{S^{n-1}} u\,u_q\,d\sigma\right) \\
&= E(u) + t_0^2 E(u_q) - 2t_0^2 E(u_q) = E(u) - t_0^2 E(u_q).
\end{align*}
Since $E(u) = E(u_q) = \gamma_n|S^{n-1}|$, it follows that
\[
\|v\|^2 = (1-t_0^2)\gamma_n|S^{n-1}|, \qquad t_0^2 = 1 - \frac{\|v\|^2}{\gamma_n|S^{n-1}|}.
\]
Expanding $\bar{J}_\tau(u)$ in powers of $v$:
\begin{align*}
\bar{J}_\tau(u) \leq{}& t_0^{\tau+1}\int_{S^{n-1}}\bar{H}\,u_q^{\tau+1}\,d\sigma + (\tau+1)\,t_0^{\tau}\int_{S^{n-1}}\bar{H}\,u_q^\tau\,v\,d\sigma\\
&+ \frac{\tau(\tau+1)}{2}\,t_0^{\tau-1}\int_{S^{n-1}}\bar{H}\,u_q^{\tau-1}\,v^2\,d\sigma + o(\|v\|^2).
\end{align*}
The remainder is $o(\|v\|^2)$, and not merely $O(\|v\|^2)$, because the first neglected order is $\tau+1>2$; this matters, as the computation at the end of this step shows. Since $t_0^2 = 1-\|v\|^2/(\gamma_n|S^{n-1}|)$, the factors $t_0^{\tau}$ and $t_0^{\tau-1}$ equal $1+O(\|v\|^2)$ and may be absorbed into that remainder; we drop them below.
Using the value of $t_0$ computed above, the first term on the right-hand side is bounded above by
\begin{equation}\label{eq:term1}
\begin{aligned}
t_0^{\tau+1}\int_{S^{n-1}}\bar{H}\,u_q^{\tau+1}\,d\sigma
&\leq \left(1-\frac{\tau+1}{2}\frac{\|v\|^2}{\gamma_n|S^{n-1}|}\right)H(0)|S^{n-1}|\bigl(1-k_1|\tilde{q}|^\alpha - k_1\lambda^\alpha\bigr)\\
&\quad + o(\|v\|^2).
\end{aligned}
\end{equation}
By the orthogonality between $v$ and $u_q^\tau$ (\cref{lem:orto}) and \cref{lem:cmasa},
\begin{equation}\label{eq:term2}
\int_{S^{n-1}}\bar{H}\,u_q^\tau\,v\,d\sigma = \int_{B_{2\rho_0}(0)}(\bar{H}-m_0)\,u_q^\tau\,v\,d\sigma \leq C_4\|v\|(|\tilde{q}|^\alpha+\lambda^\alpha).
\end{equation}
For the quadratic term, we use that $T_\varphi v$ is $E$-orthogonal to constants and coordinate functions (\cref{lem:orto2}). Recall that the Steklov eigenvalues of the ball, $\Delta\psi = 0$ in $B^n$ with $\partial\psi/\partial\eta = \sigma\psi$ on $S^{n-1}$, are $\sigma_k = k$, the eigenfunctions being the harmonic polynomials homogeneous of degree $k$ about the center $(0,\ldots,0,1)$: the constants ($\sigma=0$), the centered coordinate functions $x_1,\ldots,x_{n-1},x_n-1$ ($\sigma=1$), and then $\sigma=2$. They have the variational characterization
\[
\sigma_{k+1} = \inf\left\{\frac{\int_{B^n}|\nabla u|^2\,dv}{\int_{S^{n-1}} u^2\,d\sigma} \;:\; u\in H^1(B^n)\setminus\{0\},\ \int_{S^{n-1}} u\,\psi_j\,d\sigma = 0 \ \text{ for } j\leq k\right\},
\]
the orthogonality being in $L^2(S^{n-1})$. For a Steklov eigenfunction $\psi$ with eigenvalue $\sigma$, integration by parts gives $\langle u,\psi\rangle_E = (\sigma+\gamma_n)\int_{S^{n-1}}u\,\psi\,d\sigma$, so $E$-orthogonality to $\psi$ is equivalent to $L^2(S^{n-1})$-orthogonality. Since $x_n-1$ lies in the span of $1$ and $x_n$, the conclusion of \cref{lem:orto2} is $L^2(S^{n-1})$-orthogonality to the full $\sigma=0$ and $\sigma=1$ eigenspaces, so $T_\varphi v$ competes in the characterization of $\sigma = 2$. Therefore there exists $c>0$ (indeed $c=1$) such that
\[
1 + c \leq \frac{\int_{B^n}|\nabla T_\varphi v|^2\,dv}{\int_{S^{n-1}}(T_\varphi v)^2\,d\sigma}.
\]
Adding $\gamma_n$ to both sides and using $E(T_\varphi v) = E(v)$:
\[
\|v\|^2 = \|T_\varphi v\|^2 = E(T_\varphi v) \geq (\gamma_n+1+c)\int_{S^{n-1}}(T_\varphi v)^2\,d\sigma.
\]
To relate the quadratic form to $v$, we use the conformal covariance of $T_\varphi$: by the change of variables $x\mapsto\varphi(x)$ and the definition of $T_\varphi$,
\begin{align*}
\int_{S^{n-1}} u_q^{\tau-1}\,v^2\,d\sigma &= \int_{S^{n-1}} u_q^{\tau-1}(\varphi(x))\,v^2(\varphi(x))\,\det(d\varphi(x))\,d\sigma \\
&= \int_{S^{n-1}} \frac{(T_\varphi u_q)^{\tau-1}}{[\det(d\varphi)]^{\frac{n-2}{2(n-1)}(\tau-1)}}\;\frac{(T_\varphi v)^2}{[\det(d\varphi)]^{\frac{n-2}{2(n-1)}\cdot 2}}\;\det(d\varphi)\,d\sigma \\
&= \int_{S^{n-1}} (T_\varphi v)^2\,d\sigma,
\end{align*}
where the exponents cancel since $\frac{n-2}{2(n-1)}(\tau-1+2) = 1$. The inequality $\bar{H}\leq H(0)$ holds by construction for $\rho_0$ small, since $H(r)=H(0)-a r^{\alpha}+k(r)$ near $\mathcal{O}$ and $m_0=H|_{\partial B_{2\rho_0}(0)}<H(0)$. Therefore
\begin{equation}\label{eq:term3}
\int_{S^{n-1}}\bar{H}\,u_q^{\tau-1}v^2\,d\sigma \leq H(0)\int_{S^{n-1}}(T_\varphi v)^2\,d\sigma \leq \frac{H(0)}{\gamma_n+1+c}\|v\|^2.
\end{equation}
Combining~\eqref{eq:term1} and~\eqref{eq:term3}, the coefficient of $\|v\|^2$ is
\[
H(0)\Bigl[-\frac{\tau+1}{2\gamma_n}+\frac{\tau(\tau+1)}{2(\gamma_n+1+c)}\Bigr]
= -\frac{(\tau+1)\,c\,H(0)}{2\gamma_n(\gamma_n+1+c)}.
\]
The identity behind this is $\tau\gamma_n = \frac{n}{n-2}\cdot\frac{n-2}{2} = \frac{n}{2} = \gamma_n+1$, so that $\frac{\tau}{\gamma_n+1} = \frac{1}{\gamma_n}$ \emph{exactly}: at the level of the first nonzero Steklov eigenvalue the quadratic term is in exact balance with the concentration term, and the whole negative margin comes from the spectral gap $c$. With $c=1$ the coefficient is $-\frac{(\tau+1)H(0)}{2\gamma_n(\gamma_n+2)}$, a fixed negative constant, which is why the remainder in the expansion above had to be $o(\|v\|^2)$. The cross term~\eqref{eq:term2} is absorbed by Young's inequality,
\[
C_4\|v\|\bigl(|\tilde{q}|^{\alpha}+\lambda^{\alpha}\bigr) \leq \varepsilon\|v\|^2 + C(\varepsilon)\bigl(|\tilde{q}|^{2\alpha}+\lambda^{2\alpha}\bigr),
\]
with $\varepsilon$ small enough not to consume the coefficient just computed; the resulting terms are of higher order in $|\tilde{q}|$ and $\lambda$, because $2\alpha>\alpha$. Hence there exists $b_0>0$ such that
\begin{equation}\label{eq:Jbar-final}
\bar{J}_\tau(u) \leq H(0)|S^{n-1}|\bigl[1-b_0(|\tilde{q}|^\alpha+\lambda^\alpha+\|v\|^2)\bigr].
\end{equation}
Now combining Steps~1 and~2:
\begin{align*}
J_p(u) &\leq |\bar{J}_p(u)-\bar{J}_\tau(u)| + |J_p(u)-\bar{J}_p(u)| + \bar{J}_\tau(u) \\
&\leq o_p(1) + C_3\lambda^{(n-1)-\delta_p} + C\|v\|^{p+1}\\
&\qquad + H(0)|S^{n-1}|\bigl[1-b_0(|\tilde{q}|^\alpha+\lambda^\alpha+\|v\|^2)\bigr].
\end{align*}
Since $u\in\partial\Sigma$ forces at least one of $|q|$, $\|v\|$ to equal $\rho_0$, letting $p\to\tau$ gives $J_p(u)\leq H(0)|S^{n-1}|-\delta_0$ for $p$ sufficiently close to $\tau$.
\end{proof}

\subsection{Existence for the subcritical problem}\label{sec:existence}

Fix $1<p<\tau$. Since $p+1<\frac{2(n-1)}{n-2}$, the trace embedding $H^1(B^n)\hookrightarrow L^{p+1}(S^{n-1})$ is compact; consequently $J_p$ is of class $C^1$ on $H^1(B^n)$ and weakly continuous on bounded sets, and $S$ is a $C^1$ hypersurface, $E$ being a positive definite quadratic form. From this, $J_p$ satisfies the Palais--Smale condition on $S$ at every positive level: a sequence $\{u_m\}\subset S$ along which $J_p$ is bounded and the constrained differential tends to zero is bounded in $H^1(B^n)$ by the constraint itself, has bounded Lagrange multipliers by the computation in \cref{lem:multiplier}, and converges strongly along a subsequence because the trace embedding is compact. This is what makes the deformation lemma, and with it the mountain pass theorem, available on $S$; it is also the only role played by the restriction $p<\tau$ in this section. We first record how a constrained critical point produces a solution of~\eqref{eq:subcrit}.

\begin{lem}\label{lem:multiplier}
If $u\in S$ is a critical point of $J_p$ on $S$ with $J_p(u)>0$, then $w = \Lambda u$ solves~\eqref{eq:subcrit} with
\[
\Lambda = \Bigl(\frac{|S^{n-1}|}{J_p(u)}\Bigr)^{1/(p-1)}.
\]
In particular $u\geq 0$, and two-sided bounds $0<c_*\leq J_p(u)\leq C^*$ make $\Lambda$ bounded and bounded away from zero uniformly as $p\to\tau$.
\end{lem}

\begin{proof}
The Lagrange condition reads $(p+1)\int_{S^{n-1}}H (u_+)^{p}\varphi\,d\sigma = 2\mu\bigl(\int_{B^n}\nabla u\cdot\nabla\varphi\,dv+\gamma_n\int_{S^{n-1}}u\varphi\,d\sigma\bigr)$ for all $\varphi\in H^1(B^n)$. Testing with $\varphi = u$ gives $(p+1)J_p(u) = 2\mu E(u) = 2\mu\gamma_n|S^{n-1}|$, so $\mu>0$ exactly when $J_p(u)>0$. Testing with $\varphi = u_- = \max\{-u,0\}$: the left side vanishes because $(u_+)^p\,u_-\equiv 0$, while $\int\nabla u\cdot\nabla u_-\,dv+\gamma_n\int u\,u_-\,d\sigma = -E(u_-)$, so $\mu\,E(u_-)=0$ and $u\geq 0$. Hence $u$ weakly solves $\Delta u = 0$ with $\partial u/\partial\eta+\gamma_n u = \gamma_n\bigl(|S^{n-1}|/J_p(u)\bigr)H u^{p}$ on $S^{n-1}$, and scaling by $\Lambda$ as stated absorbs the constant.
\end{proof}

Since $H'(r)$ changes sign in the region where $H>0$, the function $H$ has either exactly one positive local maximum or at least two, and we treat the two cases separately. The configuration in the first case is constrained as follows.

\begin{lem}\label{lem:config}
Suppose $H'$ changes sign in $\{H>0\}$ and $H$ has exactly one positive local maximum $r_1$. Then $r_1\in(0,\pi)$, and each pole is either a point where $H\leq 0$ or a positive local minimum of $H$.
\end{lem}

\begin{proof}
If there were $s_1<s_2$ in $\{H>0\}$ with $H'(s_1)<0<H'(s_2)$, then $H$ would attain a local maximum with positive value at some point of $[0,s_1)$ and another at some point of $(s_2,\pi]$, contradicting uniqueness. Hence there are $s_1<s_2$ in $\{H>0\}$ with $H'(s_1)>0>H'(s_2)$, and the maximum of $H$ over $[s_1,s_2]$ is attained at an interior point with positive value; this point is $r_1$, so $r_1\in(0,\pi)$. A pole with positive value is a critical point of $H$ in $\{H>0\}$, so~\eqref{eq:flat} holds there with $a\neq 0$: if $a<0$ the pole would be a second positive local maximum, which is excluded; hence $a>0$ and the pole is a positive local minimum.
\end{proof}

\medskip
\noindent\textbf{Case (i): exactly one positive local maximum.}
Let
\[
S_r = \{u\in S : u\circ g = u \ \text{for all } g\in SO(n-1)\},
\]
where $SO(n-1)$ acts on $B^n$ by rotations about the $x_n$-axis, and fix a rotationally symmetric $\phi_0\in S_r$ supported in a band of latitudes where $H\geq c>0$, so that $J_p(\phi_0)\geq c_*>0$ uniformly for $p$ near $\tau$. The maximum of $J_p$ over $S_r$ is attained. Indeed, the fixed-point set of the action is a closed subspace of $H^1(B^n)$, hence weakly closed, and so is its intersection $B_r$ with the ball $\{E\leq\gamma_n|S^{n-1}|\}$; since $J_p$ is weakly continuous on bounded sets, it attains its maximum over $B_r$, and that maximum is at least $J_p(\phi_0)>0$. By the homogeneity $J_p(tu)=t^{p+1}J_p(u)$, a maximizer $u_p$ must satisfy $E(u_p)=\gamma_n|S^{n-1}|$, since otherwise rescaling by some $t>1$ would keep it in $B_r$ while strictly increasing $J_p$. Hence $u_p\in S_r$ and
\begin{equation}\label{eq:cri}
c_p^{\,r} := \max_{S_r} J_p \;=\; J_p(u_p) \;\geq\; J_p(\phi_0) \;\geq\; c_*>0.
\end{equation}
The maximization is carried out on the ball first because $S$, being a sphere in $H^1(B^n)$, is not weakly closed. The group $SO(n-1)$ acts on $H^1(B^n)$ by linear isometries leaving $E$ and $J_p$ invariant, so by the principle of symmetric criticality~\cite{palais} the constrained maximizer $u_p$ is a critical point of $J_p$ on all of $S$. In the present setting this needs no appeal to the general principle: writing $w = \nabla J_p(u_p)-\mu\nabla E(u_p)$ with $\mu$ chosen so that $\langle w,u_p\rangle_E=0$, invariance of $J_p$ and $E$ under an action by $E$-isometries places $w$ in the fixed-point subspace, while maximality over $S_r$ makes $w$ orthogonal to that subspace; hence $\langle w,w\rangle_E = 0$ and $w=0$. By \cref{lem:multiplier} and~\eqref{eq:cri}, a positive multiple $\Lambda(p)u_p$ solves~\eqref{eq:subcrit}; it is rotationally symmetric, nonnegative, not identically zero, and hence positive by the strong maximum principle and Hopf's boundary point lemma, as in the proof of \cref{thm:main} below.

\medskip
\noindent\textbf{Case (ii): at least two positive local maxima.}
Let $r_1$ and $r_2$ be the two positive local maxima with the \emph{least values} of $H$. By \cref{prop:sup,prop:bdy}, there exist two disjoint open sets $\mathring{\Sigma}_1$, $\mathring{\Sigma}_2\subset S$, functions $\psi_i\in\mathring{\Sigma}_i$, a number $p_0<\tau$, and $\delta>0$ such that for all $p\geq p_0$:
\[
J_p(\psi_i) > H(r_i)|S^{n-1}|-\frac{\delta}{2}, \quad i=1,2,
\]
and
\begin{equation}\label{eq:bdybound}
J_p(u) \leq H(r_i)|S^{n-1}|-\delta \quad \text{for all } u\in\partial\Sigma_i, \quad i=1,2.
\end{equation}

Define the path family $\Gamma = \{\gamma\in C([0,1],S) : \gamma(0)=\psi_1,\;\gamma(1)=\psi_2\}$ and the minimax value
\[
c_p = \sup_{\gamma\in\Gamma}\;\min_{u\in\gamma} J_p(u).
\]
By the mountain pass theorem, there exists a critical point $u_p$ of $J_p$ on $S$ with $J_p(u_p)=c_p$. Every path in $\Gamma$ crosses both $\partial\Sigma_1$ and $\partial\Sigma_2$, so by~\eqref{eq:bdybound} and the choice of $r_1$, $r_2$ as the maxima of least value,
\begin{equation}\label{eq:cpbound}
J_p(u_p) \leq \min_k\bigl\{H(r_k)|S^{n-1}|-\delta\bigr\},
\end{equation}
where the minimum is taken over all positive local maxima $r_k$ of $H$. For the lower bound we use a two-bubble path.

\begin{lem}[Two-bubble path]\label{lem:lower}
Let $r_i\neq r_j$ be positive local maxima of $H$. For every $\varepsilon>0$ there exist $\lambda_0>0$ and $p_0<\tau$ such that for $p_0\leq p<\tau$ there is a path $\gamma_0\in C([0,1],S)$ from $u_{\lambda_0,r_i}$ to $u_{\lambda_0,r_j}$ with
\[
\min_{u\in\gamma_0} J_p(u) \;\geq\; 2^{(1-p)/2}\,\min\{H(r_i),H(r_j)\}\,|S^{n-1}|\,(1-\varepsilon)\;>\;0.
\]
Consequently $c_p \geq 2^{(1-p)/2}\min\{H(r_1),H(r_2)\}\,|S^{n-1}|(1-\varepsilon)$.
\end{lem}

\begin{proof}
Write $\psi_i = u_{\lambda_0,r_i}$, $\psi_j = u_{\lambda_0,r_j}$, and for $\theta\in[0,\pi/2]$ set $w_\theta = \cos\theta\,\psi_i+\sin\theta\,\psi_j$ and $\gamma_0(\theta) = c(\theta)\,w_\theta$, with $c(\theta)>0$ chosen so that $E(\gamma_0(\theta)) = \gamma_n|S^{n-1}|$, so $\gamma_0(\theta)\in S$; since $\psi_i,\psi_j\geq 0$, the positive part is inactive along the path. The interaction $\int_{B^n}\nabla\psi_i\cdot\nabla\psi_j\,dv+\gamma_n\int_{S^{n-1}}\psi_i\psi_j\,d\sigma = \gamma_n\int_{S^{n-1}}\psi_i^{\tau}\psi_j\,d\sigma = O(\lambda_0^{(n-2)/2})$, because $\psi_j = O(\lambda_0^{(n-2)/2})$ away from $r_j$ while $\int\psi_i^{\tau}\,d\sigma$ is bounded; hence $E(w_\theta) = \gamma_n|S^{n-1}|\bigl(1+O(\lambda_0^{(n-2)/2})\bigr)$ and $c(\theta) = 1+O(\lambda_0^{(n-2)/2})$ uniformly in $\theta$. Both bubbles are $O(\lambda_0^{(n-2)/2})$ pointwise on the set $\{H<0\}$, which is at positive distance from $r_i$ and $r_j$, so
\[
\int_{\{H< 0\}} |H|\,w_\theta^{p+1}\,d\sigma = O\bigl(\lambda_0^{\,(n-2)(p+1)/2}\bigr),
\]
while on $\{H\geq 0\}$ superadditivity of $t\mapsto t^{p+1}$ on $[0,\infty)$ gives $w_\theta^{p+1}\geq \cos^{p+1}\!\theta\,\psi_i^{p+1}+\sin^{p+1}\!\theta\,\psi_j^{p+1}$. Since $\int_{\{H<0\}}|H|\psi_i^{p+1}\,d\sigma = O(\lambda_0^{(n-2)(p+1)/2})$ as well, and $J_p(\psi_i) = H(r_i)|S^{n-1}|+o_{\lambda_0}(1)$ by the computation of \cref{prop:sup},
\[
J_p(w_\theta) \;\geq\; \bigl(\cos^{p+1}\!\theta+\sin^{p+1}\!\theta\bigr)\,\min\{H(r_i),H(r_j)\}\,|S^{n-1}|\,\bigl(1-o_{\lambda_0}(1)\bigr).
\]
Finally $\cos^{p+1}\theta+\sin^{p+1}\theta\geq 2^{(1-p)/2}$ on $[0,\pi/2]$, and $J_p(\gamma_0(\theta)) = c(\theta)^{p+1}J_p(w_\theta)$ with $c(\theta)^{p+1} = 1+O(\lambda_0^{(n-2)/2})$. Choosing $\lambda_0$ small proves the claim; the bound on $c_p$ follows since $\gamma_0\in\Gamma$ after fixing $\psi_1 = u_{\lambda_0,r_1}$, $\psi_2 = u_{\lambda_0,r_2}$.
\end{proof}

By \cref{lem:multiplier,lem:lower}, a positive scalar multiple $\Lambda(p)\,u_p$ solves~\eqref{eq:subcrit}, with $\Lambda(p)$ uniformly bounded and bounded away from $0$.

\medskip
In both cases, since $u_p\in S$ and $\Lambda(p)$ is uniformly bounded, the energy of the solutions $w_p = \Lambda(p)\,u_p$ is uniformly bounded as $p\to\tau$.

\subsection{Remarks on the hypotheses}\label{sec:hyp}

The hypotheses of \cref{thm:main} can now be discussed against the estimates that use them.

\begin{rem}\label{rem:reg}
The stated regularity is the natural one: every $\alpha\in(n-2,\,n-1)$ is non-integer, and~\eqref{eq:flat} with $a\neq 0$ forces $H'(t) = a\alpha\,t^{\alpha-1}(1+o(1))$ along the meridian through $r_0$. If $H$ were of class $C^{n-1}$ near $r_0$, Taylor's theorem applied to $H'$ would give $H'(t) = \sum_{j\leq n-2}c_j t^{j} + o(t^{n-2})$; comparison with the non-integer power $t^{\alpha-1}$ forces every $c_j$ to vanish, leaving $H' = o(t^{n-2})$ and contradicting $\alpha-1<n-2$. Thus no $C^{n-1}$ function satisfies~\eqref{eq:flat} with $a\neq 0$; in particular the flatness class contains no smooth function, which is why smoothness is not assumed. Conversely, \eqref{eq:flat}--\eqref{eq:flatplus} are compatible with $H\in C^{n-2,\gamma}$ near $r_0$, $\gamma = \alpha-(n-2)$, and the model $H(r) = H(r_0)+a|r-r_0|^{\alpha}$ realizes them. For $n=3$ the global H\"older exponent $\bar\gamma$ is required by the blow-up analysis (it supplies the modulus of continuity of $\nabla H$ away from the critical points; near them, \eqref{eq:flat} gives the modulus $\alpha-1$); every function satisfying \eqref{eq:flat} is $C^{1,\alpha-1}$ near its critical points, so the hypothesis is consistent.
\end{rem}

\newpage
\begin{rem}[The flatness window]\label{rem:window}
The two bounds in the window $n-2<\alpha<n-1$ enter at different points, and the hypothesis may be relaxed accordingly. The lower bound $\alpha>n-2$, together with~\eqref{eq:flatplus} when $n\geq 5$, is what the blow-up theory of~\cite{conformal} requires (\cref{lem:EG}); it is needed at every critical point in $\{H>0\}$ when $n\geq 4$, and not at all when $n=3$. The upper bound $\alpha<n-1$ is used only twice: at the two local maxima around which the minimax is built (\cref{jp} and \cref{prop:bdy}), and at the critical points in $\{H>0\}$ with $a(r_0)>0$, where it drives the sign argument of \cref{lem:location}. At positive local maxima not used in the minimax the upper bound is not needed, since blow-up there is excluded by the energy level~\eqref{eq:cpbound} alone. In particular, for $n=3$ the flatness condition is needed only at the two chosen maxima and at the positive critical points that are not local maxima.

The upper bound in~\eqref{eq:flat} is sharp for the estimates as we run them: the transversal integrals in \cref{jp} and in Step~3 of \cref{lem:location} converge precisely when $\alpha<n-1$, and at $\alpha=n-1$ each diverges logarithmically. The concentration term $-C\lambda^{\alpha}$ of \cref{jp} would there be replaced by $-C\lambda^{n-1}\log(1/\lambda)$, which still dominates the $O(\lambda^{n-1})$ far-field remainders, and the moments of \cref{lem:location} would acquire the same factor. We do not treat the endpoint here: absorbing the logarithm calls for a strengthened remainder condition in~\eqref{eq:flat} at such points, and the derivative count in \cref{lem:transfer} changes, since $[\alpha] = n-1$ there rather than $n-2$. Beyond the endpoint, for $\alpha>n-1$, the sign of the $\lambda^{n-1}$-order correction is no longer determined by the local coefficient $a$, and a statement of a different nature should be expected, consistent with~\cite{sharaf}.
\end{rem}

\begin{rem}[The finiteness hypothesis]\label{rem:finite}
Condition~\eqref{eq:flat} already forces the critical points in $\{H>0\}$ to be isolated, since it gives $H'(r) = a\alpha|r-r_0|^{\alpha-1}\bigl(1+o(1)\bigr)\neq 0$ for $r$ near $r_0$. Finiteness is assumed in addition because isolated critical points could otherwise accumulate on the boundary of $\{H>0\}$, with the corresponding critical values accumulating at $0$; it is used only to select the two positive local maxima of least value in Case~(ii) of \cref{sec:subcritical}, and it can be relaxed, as the next two paragraphs describe.

Let $P$ denote the set of positive local maxima of $H$. Call a pair $r\neq r'$ in $P$ \emph{admissible} if
\[
H(P)\cap\bigl[\theta_n m,\;m\bigr)=\emptyset, \qquad m=\min\{H(r),H(r')\}.
\]

\cref{thm:main} remains valid if its finiteness hypothesis is replaced by the weaker requirement that $P$ consist of a single point, or else contain an admissible pair. Indeed,~\eqref{eq:flat} gives $H'(r)=\sgn(r-r_0)|r-r_0|^{\alpha-1}(a\alpha+o(1))$ near a critical point $r_0$ in $\{H>0\}$, so such critical points are isolated; consequently $H(P)$ can accumulate only at $0$, since distinct maxima with values tending to some $v>0$ would converge to a non-isolated critical point in $\{H>0\}$. Hence $H(P)\cap[c,\infty)$ is finite for every $c>0$, so for an admissible pair $\sup\bigl(H(P)\cap[0,\theta_n m)\bigr)<\theta_n m$, and $\varepsilon$ may be chosen so small that $H(P)\cap[\theta_n m(1-2\varepsilon),\,m)=\emptyset$. Running Case~(ii) with the pair $(r,r')$, the bound~\eqref{eq:cpbound} is then replaced by the pair of bounds $c_p\leq m|S^{n-1}|-\delta$ and $c_p\geq\theta_n m|S^{n-1}|(1-2\varepsilon)$ of \cref{lem:lower}, and a blow-up point would have $H(x_0)$ in the empty set just displayed. This is the bookkeeping of \cref{thm:multi}, with hypothesis~(1) there supplied by admissibility instead of assumed.

Finiteness implies admissibility of the pair realizing $\min H(P)$, and the converse fails: $H(P)=\{4^{-j}v_0:j\geq0\}$ is infinite yet admissible, since $\theta_n\geq\frac12>\frac14$ for every $n\geq3$. The requirement cannot be dropped altogether by this route: if $H(P)$ is infinite with consecutive values in ratio $\theta_n$, every window $[\theta_n m,m)$ contains the next value down and no choice of pair yields a contradiction. We retain finiteness in \cref{thm:main} because it is checkable at a glance, because the relaxation changes the form of~\eqref{eq:cpbound} and of Case~(ii) of \cref{prop:pos} rather than only the hypothesis, and because \cref{thm:multi} presupposes it in any case.
\end{rem}

\section{A priori estimates and passage to the critical exponent}\label{sec:apriori}

Having obtained subcritical solutions $\{u_p\}$ for $p<\tau$, we now show that a subsequence converges to a solution of~\eqref{eq:main} at the critical exponent $p=\tau$. Since the $u_p$ are harmonic in $B^n$, it suffices to establish uniform bounds on $S^{n-1}$.

\begin{thm}\label{thm:bound}
Under the hypotheses of \cref{thm:main}, there exists $p_0<\tau$ such that for all $p_0<p<\tau$, the solutions of~\eqref{eq:subcrit} obtained in \cref{sec:subcritical}, in either case, are uniformly bounded.
\end{thm}

We prove the theorem by estimating the solutions on three regions; the three estimates are combined at the start of \cref{sec:completion}.

\subsection{Region I: \texorpdfstring{$H$}{H} negative and away from zero}

Following the ideas of Chen and Li~\cite{chen}, we derive a priori estimates where $H$ is negative and bounded away from zero. Via the conformal extension of the stereographic projection, the problem on the half-space $\R^n_+$ reads
\begin{equation}\label{eq:halfspace}
\begin{cases}
\Delta u = 0 & \text{in } \R^n_+, \\[4pt]
\displaystyle\frac{\partial u}{\partial\eta} = H\,u^{n/(n-2)} & \text{on } \partial\R^n_+,
\end{cases}
\end{equation}
with $u(x)\sim |x|^{2-n}$ at infinity.

\begin{prop}\label{prop:neg}
The solutions of~\eqref{eq:subcrit} are uniformly bounded in the region $\{H(x)\leq -\delta\}$ for every $\delta>0$. The bound depends on $\delta$, $\dist(\{H\leq -\delta\},\,\{H=0\})$, and $\inf H$.
\end{prop}

\begin{lem}[Harnack inequality]\label{lem:harnack}
Let $x_0\in\partial\R^n_+$ satisfy $H(x_0)<0$. Let $3\varepsilon_0 = \dist(x_0,\{H=0\})$ and suppose $H(x)\leq -\delta_0$ on $\partial' B_{2\varepsilon_0}^+(x_0)$ and $H(x)\geq -M$ on $\partial\R^n_+$. Then there exists $C=C(\varepsilon_0,\delta_0,M)$ such that for any $x_1\in\partial\R^n_+$,
\[
u(x_0) \leq C\bigl(|x_1-x_0|+\varepsilon_0\bigr)^{n-2}\,u(x_1).
\]
\end{lem}

\begin{proof}
We present the proof for $p=\tau$; for $p<\tau$ the rescaling below introduces only the bounded factor $\varepsilon_0^{(n-2)(\tau-p)/2}$ in the boundary coefficient, which changes none of the constants' dependence. We may assume $x_1\neq x_0$. Choose $\bar{x}$ on the ray from $x_1$ through $x_0$, beyond $x_0$, with $|x_0-\bar{x}|\cdot|x_1-\bar{x}|=\varepsilon_0^2$; since $|x_1-\bar{x}|\geq|x_0-\bar{x}|$, this forces $|x_0-\bar{x}|\leq\varepsilon_0$. Translate coordinates so that $\bar{x}$ is the origin and define $u_0(x) = \varepsilon_0^{(n-2)/2}\,u(\varepsilon_0 x)$, which satisfies~\eqref{eq:halfspace} with $\bar{H}(x) = H(\varepsilon_0 x)$; by hypothesis, $\bar{H}\leq-\delta_0$ on $\partial' B_1^+(0)$ and $\bar{H}\geq-M$ on $\partial\R^n_+$.

Let $v(x) = |x|^{-(n-2)}\,u_0(x/|x|^2)$ be the Kelvin transform of $u_0$. A direct calculation shows that $v$ satisfies
\begin{equation}\label{eq:kelvin}
\begin{cases}
\Delta v = 0 & \text{in } \R^n_+\setminus\{0\}, \\[4pt]
\displaystyle\frac{\partial v}{\partial\eta} = \bar{H}\!\left(\frac{x}{|x|^2}\right) v^{n/(n-2)} & \text{on } \partial\R^n_+\setminus\{0\},
\end{cases}
\end{equation}
and, since $u_0(y)\sim|y|^{2-n}$ at infinity, $v$ is bounded near the origin. Fix $A > \max\{1,(M/\delta_0)^{(n-2)/2}\}$ and set $w = A v - u_0$. On $\partial'' B_1^+$ we have $v=u_0$ (the inversion fixes the unit sphere), so $w=(A-1)u_0\geq 0$ there. We claim $w\geq 0$ in $B_1^+$. For $\sigma>0$ consider the barrier-corrected function $w_\sigma = w+\sigma(|x|^{2-n}-1)$: it is harmonic in $B_1^+\setminus\{0\}$, its Neumann data on $\partial'B_1^+\setminus\{0\}$ agree with those of $w$ (indeed $\partial_n|x|^{2-n}=(2-n)|x|^{-n}x_n$ vanishes on $\{x_n=0\}$), it coincides with $w\geq 0$ on $\partial''B_1^+$, and it tends to $+\infty$ at the origin. If $w_\sigma$ were negative somewhere, its minimum over $\overline{B_1^+}\setminus\{0\}$ would be attained at some $z_0\in\partial'B_1^+\setminus\{0\}$ (not in the interior, by the strong maximum principle; not on $\partial''$, where $w_\sigma\geq0$; not at $0$, where $w_\sigma\to+\infty$), with $w(z_0)\leq w_\sigma(z_0)<0$, that is, $0\leq A v(z_0)<u_0(z_0)$. Hopf's boundary point lemma gives $\partial w_\sigma/\partial\eta(z_0)<0$. On the other hand, from the boundary conditions,
\[
\frac{\partial w_\sigma}{\partial\eta}(z_0)
= A^{-2/(n-2)}\,\bar{H}\!\left(\frac{z_0}{|z_0|^2}\right)\bigl(A v(z_0)\bigr)^{n/(n-2)} - \bar{H}(z_0)\,u_0(z_0)^{n/(n-2)},
\]
where $-\bar{H}(z_0)u_0^{n/(n-2)}\geq\delta_0\,u_0(z_0)^{n/(n-2)}$, while the first term is nonnegative if $\bar{H}(z_0/|z_0|^2)\geq0$, and otherwise is at least $-MA^{-2/(n-2)}u_0(z_0)^{n/(n-2)}$ because $0\leq A v(z_0)<u_0(z_0)$. In either case
\[
\frac{\partial w_\sigma}{\partial\eta}(z_0) \geq \bigl(\delta_0 - MA^{-2/(n-2)}\bigr)\,u_0(z_0)^{n/(n-2)} > 0,
\]
a contradiction. Hence $w_\sigma\geq 0$ for every $\sigma>0$, and letting $\sigma\to0$, $w\geq 0$ in $B_1^+$. This implies
\[
u(x) \leq A\,\frac{\varepsilon_0^{n-2}}{|x-\bar{x}|^{n-2}}\,u\!\left(\varepsilon_0^2\,\frac{x-\bar{x}}{|x-\bar{x}|^2}+\bar{x}\right)
\]
for $x\in B_{\varepsilon_0}(\bar{x})$. By the choice of $\bar{x}$, the inversion sends $x_0$ to $x_1$, and evaluating at $x=x_0$,
\[
u(x_0)\leq A\,\varepsilon_0^{2-n}\,|x_1-\bar{x}|^{n-2}\,u(x_1)\leq A\,\varepsilon_0^{2-n}\bigl(\varepsilon_0+|x_1-x_0|\bigr)^{n-2}u(x_1),
\]
which is the stated inequality with $C=A\,\varepsilon_0^{2-n}$.
\end{proof}

\begin{proof}[Proof of \cref{prop:neg}]
Recall that the energy of $u_p$ is uniformly bounded, so $\int_{\partial' B_{\varepsilon_0}^+(x_0)} u_p^{p+1}\,d\sigma$ is uniformly bounded for any $x_0$ with $H(x_0)<0$. Applying \cref{lem:harnack}, for any $x$ in the region $\{H\leq -\delta_0\}$:
\[
u_p(x) \leq C\inf_{\bar{x}\in\partial' B_{\varepsilon_0}^+(x_0)} u_p(\bar{x}) \leq \frac{C}{|\partial' B^+_{\varepsilon_0}(x_0)|}\int_{\partial' B^+_{\varepsilon_0}(x_0)} u_p\,d\sigma.
\]
By H\"older's inequality, combining the measure factors,
\begin{align*}
\frac{C}{|\partial' B^+_{\varepsilon_0}(x_0)|}\int_{\partial' B^+_{\varepsilon_0}(x_0)} u_p\,d\sigma
&\leq \frac{C}{|\partial' B^+_{\varepsilon_0}(x_0)|^{1/(p+1)}}\left(\int_{\partial' B^+_{\varepsilon_0}(x_0)} u_p^{p+1}\,d\sigma\right)^{1/(p+1)}\\
&\leq K_2,
\end{align*}
where $K_2$ depends on $\delta_0$, $\dist(\{H\leq -\delta_0\},\{H=0\})$, and $\inf H$. Hence the solutions $\{u_p\}$ are uniformly bounded where $H$ is negative.
\end{proof}

\subsection{Region II: \texorpdfstring{$H$}{H} small and close to zero}

\begin{prop}\label{prop:small}
There exist $p_0<\tau$ and $\delta>0$ such that for all $p_0<p\leq\tau$, the solutions $\{u_p\}$ are uniformly bounded where $|H(x)|\leq\delta$.
\end{prop}

\begin{proof}
Arguing by contradiction, suppose there exist $u_i = u_{p_i}$ with $p_i\to\tau$ and points $x_i\to x_0$ with $H(x_0)=0$ and $u_i(x_i)\to\infty$. Two observations organize the proof. First, the normalization $E(u_p)=\gamma_n|S^{n-1}|$, together with the Sobolev embedding of $H^1(B^n)$, bounds the conformally natural volume integral,
\begin{equation}\label{eq:Kn}
\int_{B^n} u_p^{\,2\tau}\,dv \leq C_0
\end{equation}
uniformly in $p$, and this integral is exactly conformally invariant: if $\Phi$ is the conformal map onto $\R^n_+$ centered at $x_0$, with conformal factor $F$, and $U_i = F\cdot(u_i\circ\Phi^{-1})$, then $\int_{\R^n_+} U_i^{\,2\tau}\,dz = \int_{B^n} u_i^{\,2\tau}\,dv$, with no error term. Second, $U_i$ is harmonic in $\R^n_+$ with $\partial U_i/\partial\eta = h_i\,U_i^{\,p_i}$ on $\partial\R^n_+$, where $h_i=\gamma_n(H\circ\Phi^{-1})F_b^{\,\tau-p_i}$, $F_b$ being the boundary trace of $F$; each $h_i$ is bounded and continuous and vanishes at the origin, since $H(x_0)=0$.

Select a well-scaled point as follows. With $y_i=\Phi(x_i)\to0$, $M_i=U_i(y_i)\to\infty$, and $r_i=M_i^{-(p_i-1)/2}$, maximize $\sigma_i(y)=U_i(y)\,(r_i-|y-y_i|)^{1/(p_i-1)}$ over the closed half-ball $\overline{B^+_{r_i}(y_i)}$, say at $\xi_i$, and set $K_i=U_i(\xi_i)$, $\mu_i=K_i^{-(p_i-1)}$. Since $\sigma_i(y_i)=M_i^{1/2}\to\infty$, one finds $K_i\geq M_i\to\infty$, $\xi_i\to0$, and $U_i\leq C(n)K_i$ on half-balls about $\xi_i$ whose radii, measured in units of $\mu_i$, tend to infinity. The rescalings $v_i(x)=K_i^{-1}U_i(\mu_i x+\xi_i)$ then satisfy $v_i=1$ at the base point, $0\leq v_i\leq C(n)$, $\Delta v_i=0$, and on the rescaled boundary
\[
\frac{\partial v_i}{\partial\eta} = h_i(\mu_i x + \xi_i)\,v_i^{\,p_i},
\]
with no divergent prefactor: $\mu_i=K_i^{-(p_i-1)}$ is the exact scale at which the boundary nonlinearity is invariant. Since $\xi_i\to0$, $\mu_i\to0$, and $H(x_0)=0$, the boundary data tend to zero uniformly on compact sets. Passing to a subsequence according to whether $\dist(\xi_i,\partial\R^n_+)/\mu_i$ diverges or stays bounded, elliptic estimates and the Arzel\`a--Ascoli theorem give a limit $v_0$, bounded and harmonic on $\R^n$ in the first case, and on $\overline{\R^n_+}$ with vanishing Neumann data in the second; by Liouville's theorem (after even reflection in the second case), $v_0\equiv1$. Hence for every fixed $K$ and all large $i$, $v_i\geq\frac12$ on a ball or half-ball $Q_K$ of radius $K$ about the base point, and undoing the rescaling,
\[
\int_{\R^n_+} U_i^{\,2\tau}\,dz \ \geq\ K_i^{\,2\tau-n(p_i-1)}\int_{Q_K}v_i^{\,2\tau}\,dx \ \geq\ c_n\,2^{-2\tau}\,K^n,
\]
because $2\tau-n(p_i-1)=n(\tau-p_i)\geq0$ and $K_i\geq1$. By~\eqref{eq:Kn} and conformal invariance the left side is at most $C_0$, independently of $i$ and $K$; taking $K$ large yields a contradiction. Hence $\{u_i\}$ is uniformly bounded where $H$ is small.
\end{proof}

\subsection{Region III: \texorpdfstring{$H$}{H} positive and away from zero}

We first assemble the blow-up theory of~\cite{conformal} in the form we use. Recall from \cref{prop:neg,prop:small} that blow-up can occur only in a region $\{H\geq\delta\}$, $\delta>0$.

\newpage
\begin{lem}[Localized form of the results of~\cite{conformal}]\label{lem:transfer}
Let $\delta>0$, let $\Omega\subset S^{n-1}$ be open with $\overline{\Omega}\subset\{H\geq\delta/2\}$, and let $\{u_i\}$ be solutions of~\eqref{eq:subcrit} with $p_i\to\tau$ whose blow-up points all lie in $\Omega$. Then the following results of~\cite[Section~4]{conformal} apply to $\{u_i\}$ on $\Omega$, with $h = \gamma_n H|_{\Omega}$ in place of the globally positive smooth function assumed there.
\begin{enumerate}
\item[\textup{(T1)}] \cite[Proposition~4.11]{conformal}: for $p_i$ close to $\tau$ and $\max u_i$ large there is a finite set of points of $S^{n-1}$, each a local maximum of $u_i$, near which $u_i$ is, after the natural rescaling, $C^2$-close to a standard bubble, and away from which $u_i(x)\leq C\,\dist\bigl(x,\text{that set}\bigr)^{-1/(p_i-1)}$. No condition on $h$ is imposed beyond the standing ones.
\item[\textup{(T2)}] \cite[Theorem~4.8]{conformal}: for $n=3,4$, every isolated blow-up point is isolated simple. No condition on $h$.
\item[\textup{(T3)}] \cite[Theorem~4.10]{conformal}, in the form of the Remark following it: for $n\geq5$ the same conclusion holds provided the transferred coefficients satisfy Li's condition $(*)_{n-2}$ of~\cite[Definition~0.4]{Li}.
\item[\textup{(T4)}] \cite[Theorem~4.12]{conformal}: the set of \textup{(T1)} reduces to a single point. For $n=3$ no condition on $h$ is needed, and for $n\geq4$ it suffices that $h$ satisfy condition~$(m)$ of~\cite[Definition~4.9]{conformal} with $m>n-2$ in a neighborhood of each of its critical points.
\item[\textup{(T5)}] \cite[Proposition~4.13]{conformal}: under the hypothesis of \textup{(T4)}, if $\max u_i$ is large then $E(u_i)$ is bounded and the normalized functions lie in an arbitrarily small neighborhood of the family of standard bubbles.
\end{enumerate}
\end{lem}

\begin{proof}
Three hypotheses of~\cite[Section~4]{conformal} have to be matched.

\emph{Positivity.} The results of~\cite[Section~4]{conformal} are stated on a bounded open set $\Omega'\subset\R^n_+$ with $\partial^0\Omega'$ bounded, for coefficient functions $h_i$ on $\partial^0\Omega'$ subject to $0<A\leq h_i\leq B$; we apply them with $\Omega'$ the image of $\Omega$ under the conformal map used there. Our $H$ changes sign, but by \cref{prop:neg,prop:small} blow-up occurs only where $H\geq\delta$, and on $\Omega$ we have $\gamma_n\delta/2\leq h\leq\gamma_n\max_{S^{n-1}}H$. The standing positivity therefore holds after localizing, which is possible precisely because~\cite{conformal} works on a bounded open set rather than on all of $S^{n-1}$.

\emph{Localization.} \textup{(T1)}, \textup{(T4)} and \textup{(T5)} are stated in~\cite{conformal} for the global problem, with $h$ positive on the whole boundary sphere. Their proofs, however, act only near the points of concentration: each argument rescales around a point where the solution is large and works in half-balls about it, and the positivity of $h$ enters only through those half-balls, which by \cref{prop:neg,prop:small} lie in $\{H\geq\delta\}\subset\Omega$, where the standing bounds above hold. Away from the concentration set, \cref{prop:neg,prop:small} themselves supply the uniform bounds that close each argument, including the decay statement of \textup{(T1)} used in \cref{lem:location}. This is the same localization that Chen and Li perform on $S^n$~\cite{pree}.

\emph{Regularity.} The results are stated for smooth $h$, whereas by \cref{rem:reg} our $H$ is not smooth at the flat critical points; indeed no $C^{n-1}$ function satisfies~\eqref{eq:flat}. The arguments on~\cite[pp.~138--148]{conformal} use $h$ only through Taylor expansions of $\nabla h$ at the blow-up point, to order $s\leq[\alpha] = n-2$. Conditions~\eqref{eq:flat}--\eqref{eq:flatplus} guarantee that $\nabla^{s}H$ exists and obeys the required bounds for $s\leq n-2$ near the critical points in $\{H>0\}$, and the global hypothesis $H\in C^{n-2}(S^{n-1})$, or $H\in C^{1,\bar\gamma}(S^{n-1})$ when $n=3$, supplies them elsewhere.

\emph{Flatness.} \textup{(T4)} requires condition~$(m)$ of~\cite[Definition~4.9]{conformal} with $m>n-2$ at every critical point of $h$. Here~\eqref{eq:flat} with $a\neq 0$ gives the two-sided bounds of that definition with $m=\alpha>n-2$, and it is assumed only at the critical points lying in $\{H>0\}$; this suffices because the condition enters only through neighborhoods of the blow-up points, and those lie in $\Omega$. For $n\geq 5$ we use \textup{(T3)}: the printed form of~\cite[Theorem~4.10]{conformal} assumes the condition $(n-2)$ of that paper, whose lower bound our $H$ violates, whereas the Remark following it asserts the same conclusion under Li's one-sided condition $(*)_{n-2}$. That our $H$ satisfies $(*)_{n-2}$ is a computation: for the model term $a|s|^{\alpha}$ the required inequalities $|\nabla^{s}H|\leq L\,|\nabla H|^{(n-2-s)/(n-3)}$, $2\leq s\leq n-2$, reduce to the exponent inequality
\[
(\alpha-s)(n-3)-(\alpha-1)(n-2-s) \;=\; (s-1)\bigl(\alpha-(n-2)\bigr)\;>\;0,
\]
and~\eqref{eq:flatplus} extends them to the remainder $k$.
\end{proof}

\begin{rem}[What is quoted and what is proved]\label{rem:remark410}
For $n\geq 5$ the passage from isolated to isolated simple blow-up rests on the Remark following~\cite[Theorem~4.10]{conformal}, the proof of which is left to the reader there. What \cref{lem:transfer} establishes is that our $H$ satisfies the hypothesis $(*)_{n-2}$ of that Remark; the implication itself we quote. (Li's own Definition~0.4 labels the condition simply $(*)$; the subscript in $(*)_{n-2}$ is the notation of~\cite{conformal}.) No such issue arises for $n=3$ or $n=4$, since~\cite[Theorem~4.8]{conformal} requires no condition on $h$ at all.
\end{rem}

\begin{lem}[Escobar--Garc\'ia]\label{lem:EG}
Let $\{u_i\}\subset S$ be critical points of $J_{p_i}$ with $p_i\to\tau$ and $\max u_i\to\infty$. Under the flatness hypotheses of \cref{thm:main}, after passing to a subsequence, $\{u_i\}$ has exactly one blow-up point $x_0$, with $H(x_0)\geq\delta$; the blow-up is isolated and simple; $u_i$ concentrates as a single standard bubble $u_q$ at $x_0$ with concentration parameter $\lambda_i\to 0$; and
\begin{equation}\label{eq:level}
J_{p_i}(u_i) \longrightarrow H(x_0)\,|S^{n-1}|.
\end{equation}
\end{lem}

\begin{proof}
The rescaling argument of \cref{prop:small}, run at a point with $H(x_0)>0$, produces a bounded positive limit of the natural rescalings and hence a fixed quantum of the conformally invariant energy~\eqref{eq:Kn} at each blow-up point (\cref{prop:pos}, Step~1); so there are finitely many, all lying in $\{H\geq\delta\}$, and we may fix an open $\Omega$ with $\overline{\Omega}\subset\{H\geq\delta/2\}$ containing them all and apply \cref{lem:transfer} throughout. The structure of the blow-up set is then~\textup{(T1)}. Isolated blow-up points are isolated \emph{simple} by~\textup{(T2)} for $n=3,4$ and by~\textup{(T3)} for $n\geq 5$ (see \cref{rem:remark410}). That there is at most one blow-up point is~\textup{(T4)}: for $n=3$ no condition on $H$ is needed, and for $n\geq 4$ the required condition~$(m)$ with $m>n-2$ near each critical point is supplied by~\eqref{eq:flat} with $m=\alpha$. The single-bubble asymptotics and full concentration are~\textup{(T5)}, and~\eqref{eq:level} follows since $u_i^{p_i+1}\,d\sigma$ concentrates at $x_0$ with total mass $|S^{n-1}|+o(1)$.
\end{proof}

The next lemma is the boundary-trace analogue of the localization lemma used on $S^n$ by Chen and Li~\cite[Lemma~3.1]{pree}, where it is quoted from~\cite{Li}. Since no reference covers the trace setting, we give a proof.

\begin{lem}[Blow-up occurs only at local maxima of $H$]\label{lem:location}
In the setting of \cref{lem:EG}, the blow-up point $x_0$ is a critical point of $H$. If in addition $H$ satisfies~\eqref{eq:flat} at $x_0$ with $a=a(x_0)\neq0$, then $a<0$; that is, $x_0$ is a local maximum of $H$.
\end{lem}

\begin{proof}
Write $\eta_i = \tau-p_i\geq 0$. Let $\Phi$ be the stereographic-type conformal map of $B^n$ onto $\R^n_+$ sending $x_0$ to the origin and its antipode to infinity (the map of~\cite[p.~145]{conformal}), with conformal factor $F$ normalized by $(\Phi^{-1})^*(g_{B^n}) = F^{4/(n-2)}\,g_{\R^n_+}$, so that $U_i := F\cdot(u_i\circ\Phi^{-1})$ satisfies
\[
\Delta U_i = 0 \ \text{in } \R^n_+, \qquad \frac{\partial U_i}{\partial\eta} = h_i\,U_i^{p_i} \ \text{on } \partial\R^n_+, \qquad h_i := \gamma_n\,\widetilde{H}\,F_b^{\,\eta_i},
\]
where $\widetilde{H} = H\circ\Phi^{-1}$ and $F_b$ denotes the boundary trace of $F$\footnote{The corresponding display on~\cite[p.~145]{conformal} carries the exponent $-\delta_i$ on the conformal factor, $\delta_i=\tau-p_i$ in the notation there. Direct computation from the stated normalization of $F$ gives the exponent $+\eta_i$, as here; the discrepancy is immaterial at the critical exponent, where $\eta_i=0$.}, a smooth positive function with $\nabla' F_b(0)=0$ and $F_b(x')\sim c\,|x'|^{-(n-2)}$ at infinity. Since $u_i$ is uniformly bounded away from $x_0$ by \cref{lem:EG} and the decay estimate of~\cite[Proposition~4.11(2)]{conformal}, we have $U_i(x) = O(|x|^{2-n})$ and, by harmonicity, $\nabla U_i(x) = O(|x|^{1-n})$ as $|x|\to\infty$, uniformly in $i$.

\emph{Step 1: two global identities.} Because $\Delta U_i = 0$, for any $R>0$,
\[
0 = \int_{B_R^+}\Delta U_i\,\partial_k U_i\,dv = \int_{\partial B_R^+}\Bigl(\frac{\partial U_i}{\partial\nu}\,\partial_k U_i-\frac{\nu_k}{2}|\nabla U_i|^2\Bigr)d\sigma, \qquad 1\leq k\leq n-1,
\]
and, with the multiplier $x\cdot\nabla U_i+\gamma_n U_i$,
\[
0 = \int_{\partial B_R^+}\Bigl(\frac{\partial U_i}{\partial\nu}\bigl(x\cdot\nabla U_i+\gamma_n U_i\bigr)-\frac{x\cdot\nu}{2}|\nabla U_i|^2\Bigr)d\sigma,
\]
the interior terms cancelling exactly because $\gamma_n = \frac{n-2}{2}$. The decay above makes all contributions of the spherical part $\partial'' B_R^+$ vanish as $R\to\infty$. On the flat part, $\nu_k = 0$ for $k\leq n-1$ and $x\cdot\nu = 0$, while $\partial U_i/\partial\nu = h_i U_i^{p_i}$; integrating by parts in the tangential variables (the boundary terms at infinity vanish by the same decay) yields the exact identities
\begin{align}
\int_{\partial\R^n_+}\partial_k h_i\;U_i^{\,p_i+1}\,dx' &= 0, \qquad 1\leq k\leq n-1, \label{eq:KWtrans}\\
\Bigl(\gamma_n-\frac{n-1}{p_i+1}\Bigr)\int_{\partial\R^n_+} h_i\,U_i^{\,p_i+1}\,dx' &= \frac{1}{p_i+1}\int_{\partial\R^n_+}\bigl(x'\cdot\nabla' h_i\bigr)\,U_i^{\,p_i+1}\,dx'. \label{eq:KWdil}
\end{align}

\emph{Step 2: $x_0$ is critical.} By \cref{lem:EG}, $U_i^{\,p_i+1}\,dx'$ concentrates at $0$ with mass $c_0+o(1)$, $c_0>0$, while on $|x'|\geq\sigma$ the decay gives $\int_{|x'|\geq\sigma}(1+|x'|)\,U_i^{\,p_i+1}\,dx' = O(\lambda_i^{\,n-1})$. Since $\partial_k h_i = \gamma_n F_b^{\,\eta_i}\bigl(\partial_k\widetilde{H}+\eta_i\widetilde{H}\,\partial_k F_b/F_b\bigr)$ and $\nabla' F_b(0)=0$, identity~\eqref{eq:KWtrans} gives $\nabla'\widetilde{H}(0)\,c_0 = o(1)$, hence $\nabla H(x_0)=0$.

\emph{Step 3: the sign of $a$.} Let $Z\subset\partial\R^n_+$ be the image under $\Phi$ of the critical latitude sphere through $x_0$; it is a smooth $(n-2)$-dimensional submanifold through $0$, and~\eqref{eq:flat} gives, for $d = \dist(\cdot,Z)$,
\begin{gather*}
\widetilde{H}(x') = H(x_0)+\tilde{a}\,d(x')^{\alpha}\bigl(1+o(1)\bigr), \\
x'\cdot\nabla'\widetilde{H}(x') = \tilde{a}\,\alpha\,d(x')^{\alpha}\bigl(1+o(1)\bigr)+O\bigl(|x'|^{\alpha+1}\bigr),
\end{gather*}
near $0$, where $\tilde{a}$ has the sign of $a$; the error records the curvature of $Z$. The transversal moment computation for the concentrating bubble gives
\begin{equation}\label{eq:moment}
\int_{\partial\R^n_+\cap\{|x'|\leq\sigma\}} d(x')^{\alpha}\,U_i^{\,p_i+1}\,dx' \;\asymp\; \lambda_i^{\alpha}+D_i^{\alpha},
\end{equation}
where $D_i = \dist(y_i,Z)\to 0$ is the distance of the bubble center to $Z$; the convergence of the transversal moment is exactly the condition $\alpha<n-1$. The error moments $\int|x'|^{\alpha+1}U_i^{p_i+1}$ are $o(\lambda_i^{\alpha}+D_i^{\alpha})$, and the far field contributes $O(\lambda_i^{n-1}) = o(\lambda_i^{\alpha})$. On the left of~\eqref{eq:KWdil},
\[
\gamma_n-\frac{n-1}{p_i+1} = -\,c_n'\,\eta_i\bigl(1+O(\eta_i)\bigr), \qquad c_n' = \frac{(n-2)^2}{4(n-1)}>0,
\]
while $x'\cdot\nabla' F_b/F_b$ is bounded and vanishes at $0$, so the corresponding contribution is $o(\eta_i)$. Collecting terms,~\eqref{eq:KWdil} becomes
\[
c_1\,\tilde{a}\,\bigl(\lambda_i^{\alpha}+D_i^{\alpha}\bigr)\bigl(1+o(1)\bigr) \;=\; -\,c_2\,\eta_i\,H(x_0)\bigl(1+o(1)\bigr)+o(\eta_i)+o(\lambda_i^{\alpha}),
\]
with $c_1,c_2>0$. If $a>0$ the left side is bounded below by $c\,\lambda_i^{\alpha}>0$, while the right side is at most $o(\lambda_i^{\alpha})$ because $\eta_i\geq 0$ and $H(x_0)>0$, a contradiction for large $i$. Hence $a<0$.
\end{proof}

\begin{prop}\label{prop:pos}
There exists $p_0<\tau$ such that for all $p_0<p<\tau$ and any $\delta>0$, the solutions $\{u_p\}$ are uniformly bounded where $H(x)\geq\delta$.
\end{prop}

\begin{proof}
\emph{Step 1: blow-up structure and finiteness.} Let $\{x_i\}$ be a sequence with $u_i(x_i)\to\infty$ and $x_i\to x_0$ with $H(x_0)\geq\delta$. The selection and rescaling of the proof of \cref{prop:small} apply verbatim, since they use only the boundedness and continuity of the transferred coefficient $h_i$; this time $h_i(\mu_i x+\xi_i)\to\gamma_n H(x_0)>0$ locally uniformly. If the base points drift into the interior at the natural scale, the limit is again the constant $1$ and the volume count of \cref{prop:small} gives a contradiction; otherwise the limit $v_0$ is a bounded positive solution of
\[
\begin{cases}
\Delta v_0 = 0 & \text{in } \R^n_+, \\[4pt]
\displaystyle\frac{\partial v_0}{\partial\eta} = \gamma_n\,H(x_0)\,v_0^{n/(n-2)} & \text{on } \partial\R^n_+,
\end{cases}
\]
equal to $1$ at the base point. Since $v_0\geq\frac12$ on a fixed half-ball (by the uniform gradient bounds), each blow-up point carries a fixed quantum $c_0>0$ of the conformally invariant integral~\eqref{eq:Kn}, concentrated on shrinking, eventually disjoint neighborhoods. Hence $\{u_i\}$ has finitely many isolated blow-up points. By \cref{lem:EG,lem:location}, after a subsequence there is exactly one blow-up point $x_0$; it is simple, satisfies $H(x_0)\geq\delta$, and is a local maximum of $H$.

However, even this single blow-up is impossible.

In Case~(ii), the mountain-pass bound~\eqref{eq:cpbound} gives
\[
J_{p_i}(u_i) \leq \min_k\bigl\{H(r_k)|S^{n-1}|-\delta\bigr\}
\]
over all positive local maxima $r_k$ of $H$, while~\eqref{eq:level} gives $J_{p_i}(u_i)\to H(x_0)|S^{n-1}|$ with $x_0$ a positive local maximum, a contradiction.

In Case~(i), the solutions $u_i$ are rotationally symmetric, so their (finite) set of blow-up points is invariant under $SO(n-1)$; a blow-up point off the $x_n$-axis would carry its entire $(n-2)$-dimensional orbit, contradicting finiteness. Hence $x_0$ is a pole. Since $H(x_0)\geq\delta$, \cref{lem:config} makes $x_0$ a positive local \emph{minimum} of $H$, contradicting \cref{lem:location}.
\end{proof}

\subsection{Completion of the proof of Theorem~\ref{thm:main}}\label{sec:completion}

From \cref{prop:neg,prop:small,prop:pos}, the sequence $\{u_p\}$ is uniformly bounded on $S^{n-1}$, proving \cref{thm:bound}. Uniform bounds on $S^{n-1}$ together with elliptic regularity applied to the harmonic extension into $B^n$ yield equicontinuity of $\{u_p\}$ on $\overline{B^n}$. By the Arzel\`a--Ascoli theorem, a subsequence converges uniformly to a nonnegative solution $u$ of~\eqref{eq:main} at $p=\tau$; the limit is nontrivial because $J_\tau(u) = \lim c_p \geq c_*>0$ by \cref{lem:lower} (Case~(ii)) or~\eqref{eq:cri} (Case~(i)). It is positive in $B^n$ by the strong maximum principle, and positive on $S^{n-1}$ as well: were $u(x_0)=0$ at some $x_0\in S^{n-1}$, Hopf's boundary point lemma would give $\partial u/\partial\eta(x_0)<0$, whereas the boundary condition in~\eqref{eq:main} gives $\partial u/\partial\eta(x_0) = -\gamma_n u(x_0)+\gamma_n H(x_0)u(x_0)^{\tau}=0$. Finally $u$ is of class $C^{2,\gamma}(\overline{B^n})$ by the regularity theory of Cherrier~\cite{cherr} together with Schauder estimates, since $H$ is of class $C^1$ with H\"older continuous gradient; it is harmonic, hence smooth, in the interior. This completes the proof of \cref{thm:main}.\qed

\subsection{Multiplicity: proof of Theorem~\ref{thm:multi}}

For each adjacent pair $(r_i,r_{i+1})$, $1\leq i\leq k-1$, the constructions of \cref{sec:subcritical} apply verbatim to that pair: \cref{prop:sup,prop:bdy} provide disjoint neighborhoods $\Sigma(r_i)$, $\Sigma(r_{i+1})$, functions $\psi(r_j)\in\mathring{\Sigma}(r_j)$ and $\delta>0$ with $J_p\leq H(r_j)|S^{n-1}|-\delta$ on $\partial\Sigma(r_j)$ for $j\in\{i,i+1\}$. Let $\Gamma_i$ be the family of paths in $S$ joining $\psi(r_i)$ to $\psi(r_{i+1})$, and let $c_p^{(i)}$ be the associated minimax value, attained at a critical point $u_p^{(i)}$. Every path in $\Gamma_i$ crosses both boundaries, so
\[
c_p^{(i)} \;\leq\; m_i\,|S^{n-1}|-\delta,
\]
while \cref{lem:lower} gives, for any $\varepsilon>0$ and $p$ close to $\tau$,
\[
c_p^{(i)} \;\geq\; 2^{(1-p)/2}\,m_i\,|S^{n-1}|(1-\varepsilon)\;\geq\;\theta_n\,m_i\,|S^{n-1}|(1-2\varepsilon)\;>\;0.
\]
In particular \cref{lem:multiplier} produces solutions $\Lambda(p)u_p^{(i)}$ of~\eqref{eq:subcrit} with uniformly bounded energy.

Suppose $\{u_p^{(i)}\}$ blows up along some $p_l\to\tau$. By \cref{lem:EG,lem:location}, the blow-up point is a positive local maximum $r_j$ of $H$ and $J_{p_l}(u_{p_l}^{(i)})\to H(r_j)|S^{n-1}|$. The two displays confine this limit to the interval $\bigl[\theta_n m_i(1-2\varepsilon)|S^{n-1}|,\; m_i|S^{n-1}|-\delta\bigr]$. For $j\in\{i,i+1\}$ we have $H(r_j)\geq m_i$, above the interval; for $j\notin\{i,i+1\}$, hypothesis~(1) of \cref{thm:multi} keeps $H(r_j)$ outside $[\theta_n m_i,\,m_i)$, hence outside the interval once $\varepsilon$ is small, since there are finitely many maxima. This contradiction shows $\{u_p^{(i)}\}$ is uniformly bounded in Region~III; combined with \cref{prop:neg,prop:small} and the argument of the preceding subsection, a subsequence converges to a positive solution $u^{(i)}$ of~\eqref{eq:main} with
\[
J_\tau\bigl(u^{(i)}\bigr) = \lim_{p\to\tau} c_p^{(i)} \;\in\; \bigl[\theta_n m_i(1-2\varepsilon)|S^{n-1}|,\; m_i|S^{n-1}|-\delta\bigr].
\]
By hypothesis~(2), these $k-1$ intervals are pairwise disjoint for $\varepsilon$ small, so the critical values $J_\tau(u^{(i)})$ are pairwise distinct and the solutions $u^{(1)},\ldots,u^{(k-1)}$ are pairwise distinct.\qed

\section*{Declarations}

\noindent\textbf{Funding.} The authors did not receive support from any organization for the submitted work.

\noindent\textbf{Competing interests.} The authors have no competing interests to declare that are relevant to the content of this article.

\noindent\textbf{Author contributions.} \'A.~Ortiz conceived the study, developed the results and their proofs, and wrote the manuscript. G.~Garc\'ia contributed the analytic framework on which the argument builds and reviewed the manuscript. Both authors read and approved the final version.

\noindent\textbf{Data availability.} Data sharing is not applicable to this article, as no datasets were generated or analysed.


\end{document}